\documentclass{article}	
\usepackage[utf8]{inputenc}
\usepackage{appendix}
\usepackage{graphics,epsfig,graphicx,multicol,psfrag}
\usepackage{caption}
\usepackage{subcaption}
\usepackage[usenames,dvipsnames]{pstricks}
\usepackage{epsfig}
\usepackage{pst-grad} 
\usepackage{pst-plot} 
\usepackage[space]{grffile} 
\usepackage{etoolbox} 
\makeatletter 
\patchcmd\Gread@eps{\@inputcheck#1 }{\@inputcheck"#1"\relax}{}{}
\makeatother
\usepackage{comment}
\usepackage{mathtools}
\usepackage{amsmath}
\usepackage{amssymb}
\usepackage{amsthm}
\usepackage{ifthen}
\usepackage{color}
\usepackage[margin=1.4in]{geometry}
\usepackage[UKenglish]{babel}
\usepackage{polski}
\usepackage{appendix}

\allowdisplaybreaks

\usepackage[shortlabels]{enumitem}

\newtheorem{theorem}{Theorem}[section]
\newtheorem{definition}{Definition}[section]
\newtheorem{lemma}{Lemma}[section]

\newtheorem{proposition}{Proposition}[section]
\newtheorem{remark}{Remark}[section]

\numberwithin{equation}{section}
\numberwithin{theorem}{section}
\numberwithin{proposition}{section}
\numberwithin{lemma}{section}
\numberwithin{remark}{section}
\newcommand{\intav}[1]{\mathchoice {\mathop{\vrule width 6pt height 3 pt depth  -2.5pt
\kern -8pt \intop}\nolimits_{\kern -6pt#1}} {\mathop{\vrule width
5pt height 3  pt depth -2.6pt \kern -6pt \intop}\nolimits_{#1}}
{\mathop{\vrule width 5pt height 3 pt depth -2.6pt \kern -6pt
\intop}\nolimits_{#1}} {\mathop{\vrule width 5pt height 3 pt depth
-2.6pt \kern -6pt \intop}\nolimits_{#1}}}

\def\polhk#1{\setbox0=\hbox{#1}{\ooalign{\hidewidth\lower1.5ex\hbox{`}\hidewidth\crcr\unhbox0}}}

\newcommand{\dvg}{\operatorname{div}}

\newcommand{\osc}{\operatorname{osc}}

\newcommand{\loc}{\operatorname{loc}}

\newcommand{\dsty}{\displaystyle}

\newcommand{\os}{\vbox{\hrule \hbox{\vrule 
height.6em depth0pt 
\hskip.6em \vrule height.6em depth0em}
\hrule}} 

\newcommand{\df}[1]{\buildrel\mbox{\small def}\over{#1}}

 \newcommand{\Rr}{\mathbb R}

\newcommand{\be}{\beta}

\newcommand{\Ga}{\Gamma}

\newcommand{\vp}{\varphi}
\newcommand{\om}{\omega}
\newcommand{\Om}{\Omega}

\newcommand{\lm}{\lambda}

\newcommand{\Nn}{\mathbb N}

\renewcommand{\div}{\operatorname{div}}

\newcommand{\dist}{\operatorname{dist}}

\let\TeXchi\chi
\newbox\chibox
\setbox0 \hbox{\mathsurround0pt $\TeXchi$}
\setbox\chibox \hbox{\raise\dp0 \box 0 }
\def\chi{\copy\chibox}

\newtheoremstyle{break}
  {\topsep}{\topsep}%
  {\itshape}{}%
  {\bfseries}{}%
  {\newline}{}%
\theoremstyle{break}

\theoremstyle{definition}

\newcommand{\rn}{\mathbb{R}^N}

\usepackage{setspace}
\DeclareMathOperator*{\esssup}{ess\,sup}
\DeclareMathOperator*{\essinf}{ess\,inf}
\DeclareMathOperator*{\essosc}{ess\,osc}

\title{Boundary Estimates for Doubly Nonlinear\\ Parabolic Equations}

\author{Ugo Gianazza\\
Dipartimento di Matematica ``F. Casorati", 
Universit\`a di Pavia\\ 
via Ferrata 5, 27100 Pavia, Italy\\
email: {\tt gianazza@imati.cnr.it}\\
\and
David Jesus\\
Dipartimento di Matematica,
Universit\`a di Bologna\\ 
Piazza di Porta San Donato 5, 40126 Bologna, Italy\\
email: {\tt david.jesus2@unibo.it}
}

\date{}

\begin{document}

\maketitle

\begin{abstract}
We consider non-negative, weak solutions to the doubly nonlinear parabolic equation 
\[
\partial_t u^q-\dvg(|Du|^{p-2}Du)=0
\]
in the super-critical fast diffusion regime $\dsty 0<p-1<q<\frac{N(p-1)}{(N-p)_+}$. We show that when solutions vanish continuously at the Lipschitz boundary of a parabolic cylinder $\Omega_T$, they satisfy proper Carleson  estimates. Assuming further regularity for the boundary of the domain $\Omega_T$, we obtain a power-like decay at the boundary and a boundary Harnack inequality.
\end{abstract}
\section{Introduction}
We are interested in the boundary behavior of non-negative solutions to the following doubly nonlinear parabolic equation
\begin{align}\label{eq_doubly}
    \partial_t(u^q)-\dvg(|Du|^{p-2} Du)=0, \quad\mbox{ weakly in }\, \Omega_T,
\end{align}
where $0<p-1<q<\frac{N(p-1)}{(N-p)_+}$, $\Om\subset\rn$ is an open, bounded domain, whose regularity is be specified in the following, and $\Om_T=\Om\times(0,T]$ for some $T>0$.

\subsection{Statement of the Result}
Our first result Theorem~\ref{Thm:Carleson} is a Carleson inequality for non-negative solutions to \eqref{eq_doubly} which vanish continuously at the boundary. Here, it suffices to assume that the boundary $\partial\Om$ is Lipschitz. The proof requires a flattening of the boundary, which changes our prototype equation into a more general quasi-linear parabolic one with the same structure; thus, we devote the Appendix to the study of this equation, in particular, to the proof of the H\"older continuity of a particular class of signed solutions. 

Then, in order to obtain a two-sided decay at the boundary, we restrict our attention to $C^{1,1}$ boundaries and to the prototype equation, since the proof relies on the construction of barriers which work only in this case. The comparison principle obtained in \cite[Theorem~1.2]{BS} is instrumental in this step. The Carleson estimate builds upon the decay at the boundary to establish a boundary Harnack inequality, which constitutes our second main result, and is stated in Theorem~\ref{Thm:BdHarn}.

Before starting our analysis, we present some introductory definitions.

\begin{definition}\label{def-eq-global}
A non-negative measurable function $u:\Omega_T\to[0,+\infty)$ in the class
\[
u\in C(0,T;L^{q+1}(\Om))\cap L^p(0,T;W^{1,p}(\Om))
\]
is a non-negative, weak sub(super)-solution to \eqref{eq_doubly} if
\[
\iint_{\Om_T} \left[-u^q\partial_t\varphi+|Du|^{p-2}Du\cdot D\varphi\right]\,dxdt\le(\ge)0
\]
for any non-negative function
\[
\varphi\in W^{1,q+1}_o(0,T;L^{q+1}(\Om))\cap L^p(0,T;W^{1,p}_o(\Om)).
\]
A non-negative function $u$ is a non-negative, weak solution to \eqref{eq_doubly} if it is both a weak sub-solution and a weak super-solution.
\end{definition}
\vskip.2truecm
Since for the moment things are very general, we need to define in which sense a Dirichlet boundary datum is taken. Suppose we are given $\dsty g\in L^p(0,T;W^{1,p}(\Om))$, $g$ is
continuous on $\overline{\Om}_T$ with modulus of continuity
$\om_g(\cdot)$, and we require $u=g$ on $\partial\Om\times(0,T]$.
\begin{definition}\label{def-lateral-trace}
The boundary condition $u\le g$ on $\partial\Om\times(0,T]$ for a sub-solution ($u\ge g$ for a super-solution) 
is meant in the sense that $(u-g)_+(\cdot,t)\in W^{1,p}_o(\Om)$ ($(u-g)_-(\cdot,t)\in W^{1,p}_o(\Om)$) for
a.e. $t\in(0,T]$. A function $u$ which is both a weak sub-solution which satisfies $u\le g$, and
a weak super-solution which satisfies $u\ge g$, is a solution which takes the boundary datum $g$.
\end{definition}
\begin{definition}\label{def-eq-local}
A non-negative measurable function $u:\Omega_T\to[0,+\infty)$ in the class
\[
u\in C_{\loc}(0,T;L_{\loc}^{q+1}(\Om))\cap L_{\loc}^p(0,T;W_{\loc}^{1,p}(\Om))
\]
is a non-negative, local weak sub(super)-solution to \eqref{eq_doubly}, if for every $K\subset\subset\Om$ and for every sub-interval $[t_1,t_2]\subset(0,T]$ we have
\[
\int_K u^q\varphi\,dx\Big|_{t_1}^{t_2}+\iint_{K\times(t_1,t_2)} \left[-u^q\partial_t\varphi+|Du|^{p-2}Du\cdot D\varphi\right]\,dxdt\le(\ge)0
\]
for any non-negative function
\[
\varphi\in W_{\loc}^{1,q+1}(0,T;L^{q+1}(K))\cap L_{\loc}^p(0,T;W^{1,p}_o(K)).
\]
A non-negative function $u$ is a non-negative, local, weak solution to \eqref{eq_doubly} if it is both a local weak sub-solution and a local, weak super-solution.
\end{definition}
\subsection{Novelty and Significance}
Parabolic equations of doubly nonlinear type arise in numerous applications of physics, for example, as a model for the flow of non-Newtonian fluids through a porous medium.
A comprehensive discussion on this subject is provided, for example, in \cite[Chapter 2]{BDGLS}, along with further references. The paper extensively investigates equation \eqref{eq_doubly}, covering aspects like the interior Harnack inequality and interior gradient estimates. Our work aims at complementing this research by addressing results about the boundary behavior, in particular, Carleson  estimates and the boundary Harnack inequality. The issue of gradient regularity up to the boundary remains, as far as we know, a very interesting open problem.

A well-known application (and motivation, at the same time) of boundary results comes from the study of free boundary problems of parabolic type, a topic that has attracted considerable attention in the scientific community. Indeed, starting from the stationary elliptic case, a Harnack-based approach to the study of the regularity of free boundaries was set up by Caffarelli in the seminal papers \cite{Caf87,Caf89} for the homogeneous two-phase case. More recently, De Silva developed a more flexible approach to the regularity of free boundaries in \cite{De11}, allowing the study of more general equations with a source term, where the main idea is an ``improvement of flatness" property for solutions, which is proved using a Harnack-type inequality. This idea sprouted into many different contexts, see for example \cite{Fer06, FS07, FS10}.

Regarding parabolic free boundaries, a pivotal challenge revolves around the intricacies of the classical Stefan problem. Following a trajectory akin to the elliptic case, Athanasopoulos, Caffarelli, and Salsa pioneered significant advancements in their seminal works \cite{ACS1, ACS2, ACS3}, using monotonicity arguments to obtain smoothness of the free boundary  (see also \cite{CS05}). Recently, De Silva, Forcilò and Savin in \cite{DFS21} managed to develop a more flexible argument, adapting the idea of De Silva in the elliptic case to allow for a non-stationary free boundary, and developed an improvement of flatness based on a Harnack inequality.

From this admittedly incomplete history of free boundary problems, it becomes evident that the boundary Harnack inequality plays a crucial role in the study of such problems. It is important to point out that another essential tool is the Hopf Lemma, which, in the case of doubly nonlinear equations, remains unknown at the moment. The reason is that the classical proof of Hopf's Lemma relies on gradient regularity up to the boundary, which has yet to be obtained.

Coming to the structure of the paper, after presenting introductory material in Section~\ref{S:Prelim}, we state and prove the Carleson estimate in Section~\ref{Sec:Carl}. Section~\ref{S:CP} is devoted to a brief discussion of the Comparison Principle, which is instrumental in the proof of the two-sided estimates considered in Section~\ref{S:Decay}. Finally, Section~\ref{S:BHI} collects the statement and the proof of the boundary Harnack inequality. The Appendices deal with the H\"older continuity up to the boundary of a particular class of signed solution to \eqref{eq_doubly}. This property plays a fundamental role in the proof of the Carleson estimate. 
\section{Preliminaries}\label{S:Prelim}
We collect some useful piece of notation and results we rely upon in the following.

We write $\mathbb{N}_0$ to denote $\mathbb{N}\cup\{0\}$. 
For a point $z_o\in \Rr^N\times \Rr$, $N\in \Nn$,  we always write $z_o=(x_o,t_o)$. By $K_\rho(x_o)$
we denote the cube in $\Rr^N$ with center $x_o\in\Rr^N$ and side length $2\rho>0$, whose faces are parallel to the coordinate planes in $\Rr^N$; analogously, by $B_\rho(x_o)$ we denote the ball in $\Rr^N$ with center $x_o\in\Rr^N$ and radius $\rho>0$. When $x_o=0$ we simply write $K_\rho$ or $B_\rho$, omitting the reference to the center. 
We use the symbol
\[
    (x_o,t_o)+Q_{R,S}:= K_R (x_o)\times (t_o-S,t_o]
\]
to denote a general ``backward'' parabolic cylinder with the indicated parameters. 
For $\theta >0$ we define ``backward'' parabolic cylinders by 
\[
    (x_o,t_o)+Q_\rho^-(\theta):= K_\rho (x_o)\times (t_o-\theta\rho^p,t_o],
\]
whereas ``forward'' cylinders are defined by
\[
    (x_o,t_o)+Q_\rho^+(\theta):= K_\rho (x_o)\times (t_o, t_o+\theta\rho^p].
\]
When $\theta =1$, we write $(x_o,t_o)+Q_\rho$. 
\subsection{Geometry of the Boundary}
All the following results are well-known, and we list them without proof, unless otherwise stated.
\begin{definition}\label{Def:Lip:Dom}
A bounded domain $\Omega\subset\rn$ is called a Lipschitz domain if there exists a small radius $r_o$ such that at every point $y\in \partial \Omega$ in an appropriate coordinate system, we have 
\begin{align*}
    \Omega\cap K_{8r_o}(y)=&\{(x',x_N)\in \Rr^N:x_N>\Phi(x')\}\cap K_{8r_o}(y),\\
    \partial \Omega\cap K_{8r_o}(y)=&\{(x',x_N)\in \Rr^N:x_N=\Phi(x')\}\cap K_{8r_o}(y),
\end{align*}
where $\Phi$ is a Lipschitz function with $\|D\Phi\|_{L^\infty}\leq L$. The constants $r_o$ and $L$ are independent of $y$.
\end{definition}

\begin{definition}\label{Def:2:1}
We say that a non-empty, bounded, open set $\Om\subset\rn$, $N\ge2$, is a domain of class $C^k$ for some $k\in\Nn\cup\{\infty\}$, provided that the following holds: for every point $x_*\in\partial\Om$, there exist $R>0$, an open interval $I\subset\Rr$ with $0\in I$, a rigid transformation $T:\rn\to\rn$ with $T(x_*)=0$, along with a function $\vp$ of class $C^k$, which maps $B_r(0)\subset\Rr^{N-1}$ into $I$ with the property $\vp(0)=0$, and such that if $\mathcal C$ denotes the open cylinder $B_R(0)\times I\subset\Rr^{N-1}\times\Rr$, then
\begin{itemize}
\item ${\mathcal C}\cap T(\Om)=\{(x',x_N)\in{\mathcal C}:\,x_N>\vp(x')\}$;
\item ${\mathcal C}\cap \partial T(\Om)=\{(x',x_N)\in{\mathcal C}:\,x_N=\vp(x')\}$;
\item ${\mathcal C}\cap (\overline{T(\Om)})^c=\{(x',x_N)\in{\mathcal C}:\,x_N<\vp(x')\}$.
\end{itemize}
\end{definition}
\begin{remark}
{\normalfont Any $C^1$ domain lies on only one side of its boundary.}
\end{remark}

\begin{theorem}
Let $\Om$ be a non-empty, bounded, open subset of $\rn$, $N\ge2$, and assume that $k\in\Nn\cup\{\infty\}$. Then $\Om$ is a $C^k$ domain if and only if for every point $x_*\in\partial\Om$ there exist an open neighborhood $U\subset\rn$ of $x_*$, a radius $r>0$, and a $C^k$ diffeomorphism $\dsty \psi=(\psi_1,\dots,\psi_N):U\to B_r(0)$ for which $\psi(x_*)=0$, and
\begin{itemize}
\item $\psi(\Om\cap U)=B_r(0)\times\rn_+$;
\item $\psi((\overline{\Om})^c\cap U)=B_r(0)\times\rn_-$;
\item $\psi(\partial\Om\cap U)=B_r(0)\times\partial\rn_+$.
\end{itemize}
\end{theorem}

\begin{definition}
Given $k\in\Nn$, a $C^k$ domain $\Om\subset\rn$ is called a domain of class $C^{k,1}$, provided the function $\vp$ appearing in Definition~\ref{Def:2:1} is of class $C^{k,1}$.
\end{definition}

\begin{lemma}
Let $\Om\subset\rn$ be a domain of class $C^k$ , and denote by $\nu$ its outward unit normal. Then $\Om$ is of class $C^{1,1}$ if and only if $\nu:\partial\Om\to\rn$ is a Lipschitz function.
\end{lemma}

\begin{definition}\label{Def:UEBC}
We say that a set $\Om\subset\rn$ satisfies a \emph{uniform exterior ball condition}, if there exists a number $r>0$ such that for any $x_*\in\partial\Om$ there is a point $x\in\rn$ for which $B_r(x)\cap \Om=\emptyset$ and $x_*\in\partial B_r(x)$. We call the supremum of all such numbers $r$ the \emph{UEBC} constant of $\Om$, and we denote it with $r_e$.
\end{definition}

\begin{definition}\label{Def:UIBC}
We say that a set $\Om\subset\rn$ satisfies a \emph{uniform interior ball condition}, if $\rn\backslash\Om$ satisfies an exterior ball condition. The \emph{UEBC} constant of $\rn\backslash\Om$ is referred to as the \emph{UIBC} constant of $\Om$, and it is denoted with $r_i$.
\end{definition}

We say that $\Om\subset\rn$ satisfies a \emph{uniform two-sided ball condition}, if $\Om$ satisfies both Definition~\ref{Def:UEBC} and Definition~\ref{Def:UIBC}. We  refer to the minimum between the UEBC and UIBC constants of $\Om$ as the uniform two-sided ball condition constant of $\Om$, and we denote it with $\bar r$.

\begin{theorem}\label{Thm:TSBC}
Any domain $\Om\subset\rn$ of class $C^{1,1}$ satisfies a uniform two-sided ball condition.
\end{theorem}

\begin{theorem}\label{Thm:BCLip}
If $\Om$ is a non-empty, bounded, open subset of $\rn$, which satisfies both a uniform exterior and a uniform interior ball conditions, then $\Om$ is a domain of class $C^{1,1}$.
\end{theorem}

\noindent For a short proof of Theorems~\ref{Thm:TSBC} and \ref{Thm:BCLip}, see, for example, \cite{Aikawa}.

\begin{theorem}\label{Thm:proj}
Let $\Om\subset\rn$ be a $C^1$ domain satisfying a UIBC with constant $r_i$. Then for every $x\in\Om$ such that $\dist(x,\partial\Om)<r_i$, there exists a unique $x_*\in\partial\Om$ with the property
\[
\dist(x,\partial\Om)=|x-x_*|.
\]
An analogous result holds for $C^1$ domains in $\rn$ satisfying the UEBC.
\end{theorem}
\subsection{H\"older Continuity of Solutions and Harnack Inequality}
In the framework of doubly nonlinear parabolic equations such as \eqref{eq_doubly}, the H\"older continuity of signed solutions is a delicate issue, which has not been completely settled yet (we refer the interested reader, for example, to \cite[Section~1.1]{BDGLS}). In the sequel, we are not going to deal with this topic in its full generality, since it goes largely beyond the aim of these notes. We limit ourselves
to consider what is relevant in order to prove our results.

In particular, the proof of the Carleson inequality stated in Theorem~\ref{Thm:Carleson} requires the qualitative knowledge of the H\"older continuity of $u$ up to the lateral boundary, and a quantitative, geometric decay of its oscillation, again up to the lateral boundary, which we state and prove in the Appendix (see Lemma~\ref{Lm:Geom-Decay}, and the final comment in Section~\ref{A:S:Geom-signed} of the Appendix).

Before coming to the Harnack inequality, we need to introduce a wider class of parabolic partial differential equations, of which \eqref{eq_doubly} is a prototype. Namely, we are interested in
\begin{align}\label{eq:structure}
    \partial_t(u^q)-\div {\mathbb A}(x,t,u,Du)=0,
\end{align}
where $\mathbb A$ satisfies the following structural conditions
\begin{align}\label{Def:structure}
    \begin{cases}
        {\mathbb A}(x,t,u,\xi)\cdot \xi \geq C_o|\xi|^p\\
        |{\mathbb A}(x,t,u,\xi)|\leq C_1|\xi|^{p-1}
    \end{cases} \mbox{ for a.e. }\,\, x\in \Omega, \forall\, u\in \Rr, \forall\, \xi\in \Rr^N,
\end{align}
where $0<C_o\le C_1$, and $0<p-1<q<\frac{N(p-1)}{(N-p)_{+}}$.

The following \emph{elliptic} form was obtained in \cite{BDGLS}.
\begin{theorem}[Harnack Inequality]\label{Thm:Harnack}
    Let $0 < p - 1 < q < N(p-1)/(N-p)_+$. There are
constants $\gamma > 1$ and $\sigma \in (0,1)$, which depend only on the data $p, q$ and $ N $, so that
if $u$ is a non-negative, continuous, weak solution to the parabolic equation \eqref{eq:structure}--\eqref{Def:structure} in $\Omega_T$ such that $u(x_o,t_o)>0$ and the following inclusion holds
\[
K_{8\rho}(x_o)\times \left(t_o-\mathcal{M}^{q+1-p}(8\rho)^p,t_o+\mathcal{M}^{q+1-p}(8\rho)^p\right)\subset \Omega_T
\]
with 
\[
\mathcal{M}:=\sup_{K_\rho(x_o)} u(\cdot,t_o),
\]
then for all
\begin{align}\label{eq_Harnack1}
    (x,t)\in K_\rho(x_o)\times \left(t_o-\sigma[u(x_o,t_o)]^{q+1-p}\rho^p, t_o+\sigma[u(x_o,t_o)]^{q+1-p}\rho^p \right)
\end{align}
we have
\begin{align}\label{eq_Harnack2}
    \gamma^{-1}u(x_o,t_o)\le u(x,t)\leq \gamma u(x_o,t_o).
\end{align}
\end{theorem}
This result owes its name to the fact that there is no waiting time, a striking difference with respect to the statements which hold for the heat equation or even for the degenerate case of the parabolic $p$-laplacian (see \cite{DGV}). It is a consequence of the fast forward and backward diffusion which dominates over the time derivative, giving the PDE a sort of  elliptic character. It is not surprising that this result is not stable under the limits $p\to 2$ and $q\to 1$. Further comments about Theorem~\ref{Thm:Harnack} can be found in \cite[Section~1.4]{BDGLS}.
\section{Carleson Estimate}\label{Sec:Carl}
In this section we deal with a \emph{Lipschitz} domain $\Om$ of constant $r_o,\,L$ (see Definition~\ref{Def:Lip:Dom}), and 
we aim at establishing a Carleson  estimate. The argument is very similar to the 
one used in the proof of \cite[Theorem 3.1]{AGS} for the $p$-parabolic case: the gist is proving a type of Harnack inequality up to the boundary which produces the desired decay, when combined with the H\"older continuity up to the same boundary. 

Before stating the main result of this section, some definitions clarify the general setting. We start by assuming that $0<u\leq M$ in $\Omega_T$. Fix $(x_o,t_o)\in S_T:=\partial \Omega\times(0,T]$. We introduce the time interval 
$I(t_o,\rho,h):=(t_o-h\rho^p,t_o+h\rho^p)$
and assume that 
\begin{align}\label{eq:incl_time}
    I(t_o,9\rho,M^{q+1-p})\subset (0,T].
\end{align}
Then, we define the set
\[
\tilde\Psi_\rho= \{|x_i-x_{o,i}|<\rho/4, \,|x_N|<2L\rho\}\times I(t_o,\rho,\eta_\rho^{q+1-p})\subset \Omega_T,
\]
here $\eta_\rho$ is the first root of
\[
    \max_{\tilde\Psi_\rho(\eta_\rho)}u=\eta_\rho.
\]
For each fixed $\rho>0$, $\eta_\rho$ is well-defined, since both functions $y_1(s)=\max_{\tilde \Psi_\rho(s)}u$ and $y_2(s)=s$ are monotone increasing and satisfy 
\begin{align*}
    \begin{cases}
        y_1(0)\geq 0,\\
        y_2(0)=0,
    \end{cases}\qquad
    \begin{cases}
        y_1(M)\leq M,\\
        y_2(M)=M.
    \end{cases}
\end{align*}
For $0<\rho<r_o$, we also define $x_\rho=(x_o', 2L\rho)$, and let $P_\rho=P_\rho(x_o,t_o)=(x_o',2L\rho,t_o)$ be such that $u(P_\rho)>0$. Note that $\dist(x_\rho, \partial \Omega)$ is of order $\rho$. The main result is the following Carleson  estimate for Lipschitz domains.
\begin{theorem}[Carleson  estimate] \label{Thm:Carleson}
    Let $u$ be a non-negative, local, bounded, weak solution to \eqref{def-eq-global} which vanishes continuously on the portion of the lateral boundary
    \[
    \left(\partial \Omega\cap \{|x_i-x_{o,i}|<2\rho, |x_N|<8L\rho\}\right)\times I(t_o,9\rho,M^{q+1-p}).
    \]
    Then, there exist constants $\gamma>0$ and $\alpha\in (0,1)$, depending only on $p,q, N, L$ such that for every $(x,t)\in \tilde \Psi_\rho$
    \begin{equation}\label{Eq:CE}
    u(x,t)\leq \gamma \left(\frac{\dist(x,\partial \Omega)}{\rho}\right)^\alpha\sup_{\tau \in I(t_o,\rho,2\eta_\rho^{q+1-p})} u(x_\rho,\tau).
    \end{equation}
\end{theorem}
\begin{remark}
\upshape More properly, \eqref{Eq:CE} should be referred to as a \emph{Carleson-type} estimate, since classical Carleson estimates are given in terms of the value of $u$ at a single point, whereas here we need to take the supremum over a proper time interval. 
\end{remark}
We proceed as in \cite{AGS}, and start by flattening the boundary. Using the notation from Definition \ref{Def:Lip:Dom}, if we define the new variables
\begin{align*}
    y_i=x_i, \, i=1, 2, ..., N-1,\quad y_N=x_N-\Phi(x'),
\end{align*}
the portion $\partial \Omega\cap \{|x_i-x_{o,i}|<2\rho, |x_N|<8L\rho\}$ coincides with a portion of the hyperplane $y_N=0$. Let $\tilde K_{2\rho}(x_o)=\{|x_i-x_{o,i}|<2\rho, |x_N|<4L\rho\}$. We orient $y_N$ so that $\Omega\cap \tilde K_{2\rho}(x_o)\subset \{y_N>0\}$. 

Since this change of variables does not affect the time variable $t$, as demonstrated in \cite{AGS} (see also \cite[Chapter X]{DiBen_1993}),  we obtain an equation like
\eqref{eq:structure}--\eqref{Def:structure} for the new variables, with two minor differences: here $\mathbb A$ does not explicitly depend on $t$, and
the two parameters $C_o, C_1$ are not given a-priori, but depend on universal constants and on $L$. 

We denote again by $x$ the new variables and say that in this case the boundary $\partial \Omega$ is \emph{flat with respect to $x_N$}. Define the following cube, which touches the boundary
\[
K^*_{2\rho}(x_o)=\{|x_i-x_{o,i}|<2\rho,\, 0<x_N<2L\rho\}.
\]
We also obtain the new set 
\begin{align*}
    \Psi_\rho^*=\Psi_\rho^*(\eta_\rho)&= \{|x_i-x_{o,i}|<2\rho, 0<x_N<4L\rho\}\times I(t_o,9\rho,\eta_\rho^{q+1-p})\subset \Omega_T,
\end{align*}
where now $\eta_\rho$ is the first root of
\[
    \max_{\Psi_\rho^*(\eta_\rho)}u=\eta_\rho.
\]
Proving Theorem \ref{Thm:Carleson} reduces to proving the following Lemma.
\begin{lemma}\label{Lm:Bry_Harn}
    Let $u$ solve \eqref{eq:structure}-\eqref{Def:structure}, $0<u\leq M$ in $\Omega_T$. Assume that $\partial \Omega$ is flat with respect to $x_N$ and that
    \begin{align}\label{eq:time_inc}
        I(t_o,9\rho,M^{q+1-p})\subset (0,T].
    \end{align}
    Suppose that $u$ vanishes continuously on 
    \[
        \partial \Omega\cap \{|x_i-x_{o,i}|<2\rho,\, x_N=0\}\times I(t_o,9\rho,M^{q+1-p}).
    \]
    There exists a constant $\gamma>0$ depending only on $p, q, N, L$ such that
    \[
    u(x,t)\leq \gamma \sup_{\tau \in I(t_o,\rho,2\eta_\rho^{q+1-p})}u(x_\rho,\tau).
    \]
\end{lemma}

\begin{proof}
    We start with a second change of variables, namely we let
    \[
    x\to \frac{x-x_o}{2L\rho}, \quad t\to \frac{1}{\eta_\rho^{q+1-p}}\frac{t-t_o}{\rho^p};
    \]
    such a change maps $K^*_{2\rho}(x_o)\times I(t_o,9\rho,2\eta_\rho^{q+1-p})$ into $\tilde Q = \{|y_i|<1/L, 0<y_N<2\}\times (-9^p,9^p]$, $x_\rho$ into $y_o=(0,...,0,1)=e_N$, $\tilde K_\rho(x_o)$ into $\tilde K_1=\{|y_i|<1/(2L), |y_N|<1\}$, and the portion $S_T\cap \overline{\Psi}^*_\rho$ of the lateral boundary into 

    \[
    \Xi=\{(y',0):|y_i|<1/L\}\times (-9^p,9^p].
    \]
    After denoting again by $(x,t)$ the transformed variables and letting $y_o=e_N$, the rescaled function
    \[
    v(x,t)=\frac{1}{\eta_\rho}u(2L\rho x+x_o,\eta_\rho^{q+1-p}\rho^pt+t_o)
    \]
    is a non-negative, bounded, weak solution to 
    \[
    \partial_t(v^q)-\div_y \tilde {\mathbb A}(x,v,D_x v)=0, \quad \mbox{ in } \tilde Q,
    \]
    where $\tilde{\mathbb A}$ still satisfies the structure conditions \eqref{Def:structure}, and $0<v\leq 1$ in $\tilde Q$. In order to simplify the proof, without loss of generality we may assume that $L=1$. 

    Since the equation has slightly changed, we denote by $\gamma_*$ and $\sigma_*$ the new constants claimed by Theorem \ref{Thm:Harnack}. We repeatedly apply \eqref{eq_Harnack1}-\eqref{eq_Harnack2}; due to all our assumptions, the only condition we need to take into account each time is that $K_{8R}(x)\subset \Omega$, where $x\in \{|x_i|<1/2, 0<x_N<1\}$, and the radius $R$ depends on the context.

    Consider the set $F_o=\{(x',x_N):|x_i|<1/2,x_N=1\}$, the points $P_{h'}$ given by
    \begin{align*}
        &x_i=h_i/8, \, h_i=-3,-2, ..., 3,\quad i=1, ..., N-1,\\
        &x_N=1,
    \end{align*}
    and the $(N-1)$-dimensional cubes $K_{1/4}(P_{h'})\cap \{x_N=1\}$. Notice that
    \begin{itemize}
        \item these $(N-1)$-dimensional cubes give a decomposition of $F_o$ of equal size;
        \item due to their size and their distance to the boundary, we have $K_2(P_{h'})\subset \Omega$ and therefore, we can always apply \eqref{eq_Harnack1}-\eqref{eq_Harnack2}.
    \end{itemize}
    Consequently, for any
    \begin{align}\label{eq:3.14}
        t_*\in \left(-\sigma_* v^{q+1-p}(y_o,0),\sigma_* v^{q+1-p}(y_o,0)\right],
    \end{align}
    by \eqref{eq_Harnack2} we have
    \[
        v(x,t_*)\leq \gamma_* v(P_{h'},t_*)
    \]
    for all $x\in K_{1/4}(P_{h'})\cap \{x_N=1\}$. On the other hand, it is easy to see that
    \[
        v(P_{h'},t_*)\leq \gamma_*^3 v(y_o,t_*),
    \]
     since we can connect $P_{h'}$ with $y_o$ using at most 3 cubes with side-length $1/4$, which guarantees the condition $K_2(P_{\bar h'})\subset \Omega$.   Therefore, for any $x\in F_o$,
    \[
    v(x,t_*)\leq \gamma_*^4 v(y_o,t_*).
    \]
    Consider the slab
    \[
    S_o=\{(x',x_N):|x_i|<1/2,\,7/8<x_N<1\}.
    \]
    As noticed above, $K_2(\bar x)\subset \Omega$ for any $\bar x\in F_o$; therefore, we can apply \eqref{eq_Harnack2} at time level $t_*$ and conclude that 
    \[
    \forall \bar x\in F_o, \forall x\in K_{1/4}(\bar x),\, v(x,t_*)\leq \gamma_* v(\bar x, t_*).
    \]
    Consequently, for all $x\in S_o$,
    \begin{align}\label{eq:3.15}
        v(x,t_*)\leq \gamma_*^5v(y_o,t_*).
    \end{align}
    Note that \eqref{eq:3.15} holds in particular for any $x\in F_1$, where
    \[
    F_1=\{(x',x_N):|x_1|<1/2, x_N=7/8\}.
    \]
    To help understand how we can iterate this argument, let us consider one more step. Let 
    \[
    x\in S_1=\{(x',x_N):|x_i|<1/2,\,(7/8)^2<x_N<7/8\}.
    \]
    There exists $\bar x\in F_1$ such that $x$ belongs to the cube $K_{7/32}(\bar x)$. Furthermore, as before, we see that $K_{7/4}(\bar x)\subset \Omega$. Thus we can apply \eqref{eq_Harnack2} at time $t_*$  and conclude that
    \[
    v(x,t_*)\leq \gamma_* v(\bar x,t_*)\leq \gamma_*^6 v(y_o,t_*).
    \]
    Therefore, we are now in a position to iterate this argument and get that
    \[
    \forall x \in S_k=\{(x',x_N):|x_i|<1/2, \,(7/8)^{k+1}<x_N<(7/8)^k\}
    \]
    there holds
    \begin{align}\label{eq:3.16}
        v(x,t_*)\leq \gamma_*^{k+5}v(y_o,t_*).
    \end{align}
    On the other hand, by \eqref{eq_Harnack1}-\eqref{eq_Harnack2}, for any $t_*$ as in \eqref{eq:3.14},
    \begin{align}\label{eq:3.17}
        v(y_o,t_*)\leq \gamma_* v(y_o,0).
    \end{align}
    Thus,
    \[
    v(x,t)\leq \gamma_*^{k+6}v(y_o,0)
    \]
    for any 
    \[
    (x,t)\in S_k\times \left[-\sigma_*v^{q+1-p}(y_o,0),\sigma_* v^{q+1-p}(y_o,0) \right].
    \]
    For any $\tau\in[-2,2]$ (and not just for $\tau=0$) we can repeat this argument verbatim, and conclude that 
    \[
    v(x,t)\leq \gamma_*^{k+6}v(y_o,\tau)
    \]
    for all $(x,t)\in S_k\times \left[-\sigma_*v^{q+1-p}(y_o,\tau),\sigma_* v^{q+1-p}(y_o,\tau) \right]$ and any $\tau \in [-2,2]$ provided $v(y_o,\tau)>0$. Hence, setting 
    \[
    \mathcal{M}_2=\sup_{\tau\in [-2,2]} v(y_o,\tau)
    \]
    yields
    \begin{align}\label{eq:3.18}
        v(x,t)\leq \gamma_*^{k+6} \mathcal{M}_2
    \end{align}
    for all $(x,t)\in S_k\times [-2,2]$. Now let 
    \begin{align*}
        &Q_*=\{|x_i|<1/4, \, 0<x_N<1, \, t\in [-1,1]\},\\
        &Q^*=\{|x_i|<1/4, \, -1<x_N<1, \, t\in [-1,1]\},\\
        &\tilde Q_*=\{|x_i|<1/2, \, 0<x_N<1, \, t\in [-2,2]\},\\
        &\tilde Q^*=\{|x_i|<1/2, \, -1<x_N<1, \, t\in [-2,2]\}.
    \end{align*}
    As it becomes clear along the proof, we have to assume $|x_i|<1/4$ in $Q_*$ for two closely related reasons:
    \begin{itemize}
        \item The first point $P_1$ we are going to choose must lie in a suitably small cylinder, so that the sequence $\{P_m\}_m$ is all contained in $\tilde Q_*$.
        \item We need to choose $k_o$ that depends only on the data.
    \end{itemize}
    We extend $v$ from $\tilde Q_*$ to $\tilde Q^*$ by odd reflection, which is still a (signed) solution (see Section~\ref{S:A:Ext} of the Appendix for the justification of this fact).

    Suppose there exists $P_1=(x_{1,1}, ..., x_{1,N},t_1)\in Q_*$ such that
    \begin{align}\label{eq:3.19}
        v(P_1)\geq \gamma_*^{k_o+6}\mathcal{M}_2.
    \end{align}

    By \eqref{eq:3.18}, we must have $0<x_{1,N}<(7/8)^{k_o+1}$, $|x_{1,i}|<1/4$, $t_1\in [-1,1]$, since $P_1\not \in S_{k}$ for any $k\leq k_o$. Consider the cylinder 
    \[
        Q(P_1)=\{|x-x_1|<2(7/8)^{k_o+1}\varepsilon^{-l},\, t_1-(2(7/8)^{k_o+1}\varepsilon^{-l})^p<t\leq t_1\},
    \]
    where $\varepsilon$ is as in Lemma~\ref{Lm:Geom-Decay}. Without loss of generality, we can assume $\varepsilon<7/8$. Now let
    \[
    \omega_o^{(1)}=\osc_{Q(P_1)}v.
    \]
    Even though we do not know the exact value of $\omega_o^{(1)}$, we can infer from the construction of $v$ by odd reflection and the normalization of $v$ that
    \begin{align}\label{eq:3.20}
        2\mathcal{M}_2\gamma_*^{k_o+6}\leq \omega_o^{(1)}\leq 2.
    \end{align}
    Provided $k_o$ is large enough, we have the inclusion $Q(P_1)\subset \tilde Q^*$. Set $\sigma_o^{(1)}=2(7/8)^{k_o+1}\varepsilon^{-l}$ and consider
    \[
    Q_o^{(1)}=K_{\sigma_o^{(1)}}(x_1)\times (t_1-\theta_o^{(1)}(\sigma_o^{(1)})^p,t_1),
    \]
    where $\theta_o^{(1)}=(\omega_o^{(1)}/A)^{q+1-p}$,
    and $A$ is the quantity claimed in Lemma~\ref{Lm:Geom-Decay}, which we can assume without loss of generality to be larger than $2$. Clearly $Q_o^{(1)}\subset Q(P_1)\subset \tilde Q^*$. By Lemma~\ref{Lm:Geom-Decay} we can build a sequence 
    \begin{equation*}
    \om_n^{(1)} = \delta\om_{n-1}^{(1)}, \quad
    \theta_n^{(1)} = \left ( \frac{\om_n^{(1)}}{A} \right )^{q+1-p}, \quad
    \sigma_n^{(1)} = \epsilon \sigma_{n-1}^{(1)},\quad\text{ and }\,\,\, Q_n^{(1)} = Q_{\sigma_n^{(1)}}^- (\theta_n^{(1)}),
    \end{equation*}
    for all non-negative integers $n$. Such a sequence satisfies
    \[
    Q_{n-1}^{(1)}\subset Q_n^{(1)}, \quad \osc_{Q_n^{(1)}}v\leq \omega_n^{(1)}.
    \]
    By iteration
    \[
    \osc_{Q_n^{(1)}}v\leq \delta^n \omega_o^{(1)}=\delta^n\osc_{Q(P_1)}v,
    \]
    which implies
    \[
    \osc_{Q(P_1)}v\geq \frac{1}{\delta^n}\osc_{Q_n^{(1)}}v.
    \]
    We choose $n=l$ and $l$ such that $\delta^{-l}>\gamma_*^5$. Then, by the choice of $\sigma_o^{(1)}$ we see that
\[
    P_1\in Q_l^{(1)}=\{|x-x_1|\leq 2(7/8)^{k_o+1},\, t_1-(2(7/8)^{k_o+1})^p<t\leq t_1\}.
\]
    Therefore, it is apparent that also $-P_1\in Q_l^{(1)}$ and thus, by the odd reflection of $v$, 
    \[
        \osc_{Q_l^{(1)}}v\geq 2 v(P_1)\geq 2\gamma_*^{k_o+6}\mathcal{M}_2,
    \]
    where we used \eqref{eq:3.19}. Hence,
    \[
    \omega_o^{(1)}\geq 2\gamma_*^{k_o+11}\mathcal{M}_2,
    \]
    which improves the previous bound \eqref{eq:3.20}. As $v$ was built by odd reflection, we conclude there exists $P_2=(x_2,t_2)=(x_{2,1}, ..., x_{2,N}, t_2)\in Q(P_1)$ such that
    \[
        v(P_2)\geq \gamma_*^{k_o+11}\mathcal{M}_2.
    \]
    As before, by \eqref{eq:3.18} we have $0<x_{2,N}<(7/8)^{k_o+6}$, $|x_{2,i}|<1/4$ and also $t_1-(2(7/8)^{k_o+1}\varepsilon^{-l})^p<t_2<t_1$. Set
    \[
    Q(P_2)=\{|x-x_2|<2(7/8)^{k_o+1}\varepsilon^{-l},\, t_2-(2(7/8)^{k_o+1}\varepsilon^{-l})^p<t<t_2\}.
    \]
    Once again, provided $k_o$ is large enough, we can assume that $Q(P_2) \subset\tilde Q^*$. By induction, we obtain a sequence $\{P_m=(x_m,t_m)\}$ such that
    \begin{align}\label{eq:3.21}
        v(P_m)\geq \gamma_*^{k_o+6+5(m-1)}\mathcal{M}_2
    \end{align}
    and
    \[
    0<x_{m,N}<\left(\frac{7}{8}\right)^{k_o+1+5(m-1)}.
    \]
    By choosing $k_o$ large enough, the sequence $\{P_m\}_m$ is all contained in the fixed cylinder $\tilde Q_*$; since $0<v<1$ and $\mathcal{M}_2\in (0,1]$, \eqref{eq:3.21} leads to a contradiction. Therefore, there exists $\tilde \gamma$ which depends only on the data, such that
    \[
    \forall (x,t)\in Q_*,\,\,\, v(x,t)\leq \tilde\gamma \mathcal{M}_2.
    \]
    We now switch back to the original variables, and we conclude that there exists a constant $\tilde \gamma$, depending only on the data, such that
    \[
    \forall (x,t)\in \{|x_i-x_{o,i}|<\rho/4,\, 0<x_N<2L\rho\}\times I(t_o,\rho,\eta_\rho^{q+1-p}),
    \]
    it holds
    \[
    u(x,t)\leq \tilde \gamma \sup_{\tau \in I(t_o,\rho,2\eta_\rho^{q+1-p})}u(x_\rho,\tau).
    \]
\end{proof}
\begin{remark}
{\normalfont A careful check of the proof shows that Theorem~\ref{Thm:Carleson} holds not only for non-negative solutions to \eqref{eq_doubly}, but also for non-negative solutions to \eqref{eq:structure}--\eqref{Def:structure}. Indeed, the proof is purely structural, and it does not rely on specific features of \eqref{eq_doubly}, which more general quasi-linear equations do not enjoy. As we are momentarily going to discuss, things change for the boundary Harnack inequality.}
\end{remark}
\section{Comparison Principle for Doubly Nonlinear Parabolic Equations}\label{S:CP}
In the following, the assumption $u\le v$ on $\partial\Om\times(0,T)$ has to be understood in the sense
\[(u-v)_+\in L^p(0,T;W^{1,p}_o(\Om));
\]
(see also Definition~\ref{def-lateral-trace}). We rely on the following result, proven in \cite[Theorem~1.2]{BS} (see also the comments by the authors on the range of parameters $q$ qnd $p$).
\begin{theorem}\label{Thm:Comp:Pr}
Let $q>0$, $p>1$ and $u$ be a non-negative, weak sub-solution and $v$ a non-negative, weak super-solution to \eqref{eq_doubly} in $\Om_T$, such that 
\[
\essinf_{\partial\Om\times(0,T)}v>0,\quad\text{ and }\quad\esssup_{\partial\Om\times(0,T)} u<\infty\,\,\text{if}\,\,q>1.
\]
If $u\le v$ on the parabolic boundary $\partial_p\Om_T$, then we have $u\le v$ a.e. in $\Om_T$.
\end{theorem}
It is worth noticing that, despite its easy and natural statement, the proof of Theorem~\ref{Thm:Comp:Pr} is definitely far from trivial.
\section{Decay at the Boundary}\label{S:Decay}
In the next two sections we need to assume more regularity on $\Om$, namely we have to deal with a \emph{$C^{1,1}$} domain. This strongly depends on the method we use, that is, on the construction of suitable barriers, which is technically feasible only for the prototype equation \eqref{eq_doubly} and for quite regular domains. It is an interesting and open problem to give a direct and structural proof of the boundary decay. This could definitely allow to consider less regular domains.

We follow the approach presented in \cite{SV_1994}, and at the same time we fix the inaccuracies the authors incurred into, when they postulated that their result holds for 
$C^{1,\alpha}$ domains with $\alpha\in(0,1)$, which unfortunately is not the case. Moreover, there are some extra technical difficulties, which are linked to the statement of the Comparison Principle of Theorem~\ref{Thm:Comp:Pr}.
Similar results were developed in \cite{adamowicz,Aikawa} in the elliptic setting, and in \cite{kuusi2014,kuusi2020} in the nonlinear parabolic one. In particular, in \cite{kuusi2014} Kuusi, Mingione and Nystr\"om use a different approach with respect to ours, and this allows them to consider more general operators of $p$-laplacian-type (although still in $C^{1,1}$ domains).

By Theorem~\ref{Thm:proj} if we take $x_o\in\Om$, and let $\Ga_{x_o}$ be the set of elements of $\partial\Om$ closest to $x_o$, i.e.
\[
\Ga_{x_o}\df=\{y\in\partial\Om:\,\,|x_o-y|=\dist(x_o,\partial\Om)\}\not=\emptyset,
\]
we have that such a $y$ is unique, provided $x_o$ is sufficiently close to the boundary.
\subsection{Control from above}\label{SS:ctrl-above}

\begin{proposition}[Control from above]\label{Prop:above}
Let $0\le u\le M$ be a weak solution to \eqref{eq_doubly} in $\Om_T$, pick $s\in(0,T)$, and fix $x_o\in\partial\Om$. Assume that $u$ vanishes on $(\partial\Om\cap B_{2\rho}(x_o))\times(s,T)$.
For every $\nu>0$ there exist $\gamma_1>0$ and $\bar s\in(0,\frac12)$, both depending on $N, q, p, \nu$, and $\|\partial \Om\|_{1,1}$, such that for all $t\in(s+\nu M^{q-p+1}\rho^p,T)$, and for all $x\in\Om\cap B_{2\bar s\rho(x_o)}$ with $\dist(x,\partial\Om)<\bar s\rho$ it holds
    \begin{align}\label{eq_upbd}
        u(x,t)\leq \gamma_1 \left(\frac{\dist(x,\partial\Om)}\rho\right) \sup_{(\Om\cap B_{2\rho}(x_o))\times(s,T)} u.
    \end{align}
\end{proposition}
\noindent Up to a possible rescaling, in the following we assume $\rho=1$. 
Let $d(x):=\dist(x,\partial \Omega)$, $\Omega_s:=\{x\in \Omega : d(x)<s\}$. Moreover, for $x\in\Om$ and $y\in \Gamma_x$, let
\begin{align*}
    y_\lambda\df=y+\lambda\frac{y-x}{|y-x|}, \quad d_\lambda=d(y_\lambda), \quad \lambda\geq 0.
\end{align*}
We have a simple geometrical result.

\begin{lemma}\label{Lm:Exterior}
    There exist $\overline{s}, \overline{\lambda}>0$ depending only upon $\Om$, such that for every $x\in \Omega_{\overline{s}}$, $y\in \Gamma_x$ and $\lambda\in (0,\overline{\lambda})$ it holds
    \begin{align*}
        \lambda=d_\lambda.
    \end{align*}
\end{lemma}
\begin{center}
\begin{figure}
\psfragscanon
\psfrag{x}{$\scriptstyle x$}
\psfrag{y}{$\scriptstyle y$}
\psfrag{ylam}{$\scriptstyle y_\lm$}
\psfrag{Omega}{$\scriptstyle \Om$}
\psfrag{lam-lame}{$\scriptstyle \lambda\le\overline{\lambda}$}
\psfrag{ball}{$\scriptstyle B_\lambda(y_\lambda)$}
\begin{center}
\includegraphics[width=.75\textwidth]{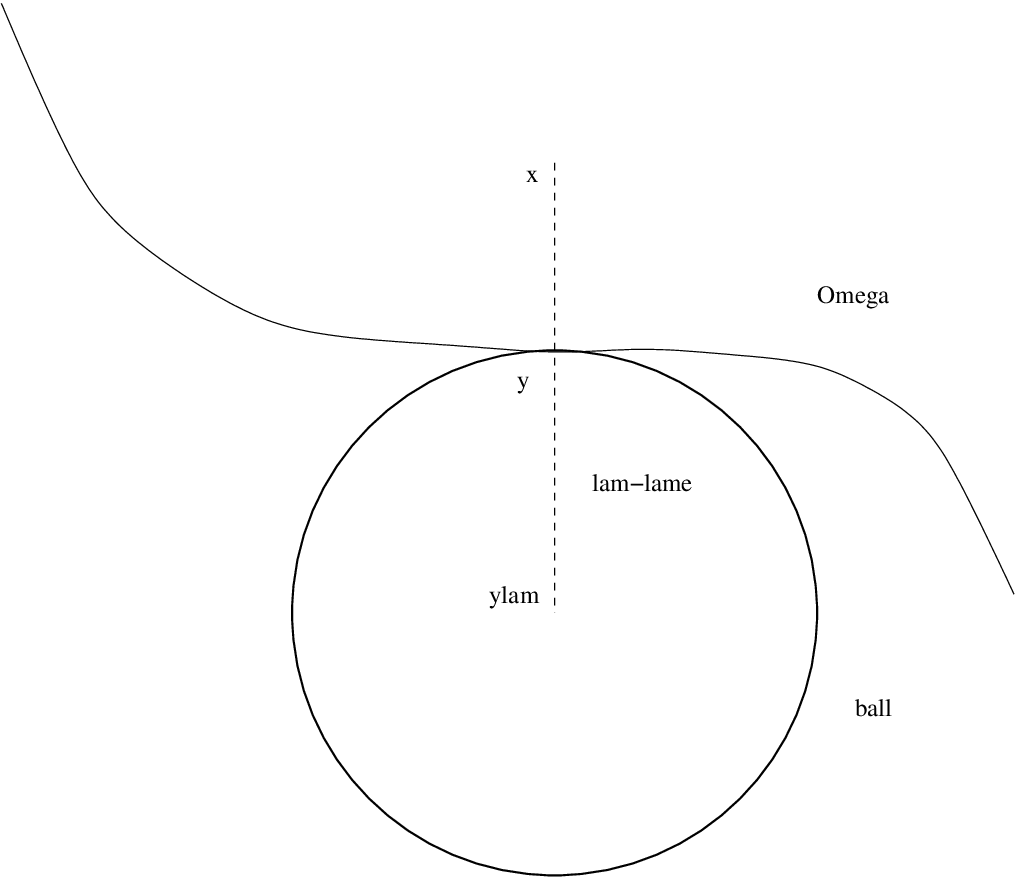}
\caption{\label{Fig:1}
}
\end{center}
\end{figure}
\end{center}
\begin{proof}
We start by taking $\overline{s}$ equal to the interior radius $r_i$, as to ensure that $\Gamma_x$ contains a single point. 
Then, by Definition~\ref{Def:UEBC}, it suffices to take $\overline{\lambda}$ equal to the exterior radius $r_e$ of the ball that touches the boundary from the exterior at $y\in \Gamma_x$; by assumption, such a radius is uniform and thus independent of $y$, and therefore also of $x$. In terms of $C^{1,1}$ regularity, if we assume we can touch the graph of $\partial \Omega$ from the exterior by paraboloids of opening $K$, then we take $\overline{\lambda}$ equal to $(\sqrt{2}-1)/K$. Figure~\ref{Fig:1} describes what happens.
\end{proof}

\noindent\textit{Proof of Proposition~\ref{Prop:above}.}
Set $\bar a=(q^{-1}-1)^+$ and fix $\nu>0$. Let $I_\nu\df=(s+\nu M^{q+1-p},T)$ and for $r\in(0,1)$ set
\begin{align*}
    M_r=\sup_{\Omega_r\times I_\nu} u(x,t)\leq M.
\end{align*}

\noindent Choose $(x_o,t_o)\in \Omega_r\times I_\nu$, with $r<r_i$ so that the corresponding $y\in \Gamma_{x_o}$ is unique. Moreover, let $I_\nu(t_o)=(t_o-\nu M^{q-p+1},t_o]\subset (s,T)$. For $c,\,k>0$ construct the barrier
\begin{equation}\label{eq_barrier}
    \begin{aligned}
        &\psi(x,t)=c M_{1/k}(1-\eta(x,t)),\\
        &\eta(x,t)=\exp\left[-k(|x-y_\lambda|-d_\lambda)\right]\exp\left[-\left(\frac{t_o-t}{M^{q-p+1}}\right)^{\bar a+1}\right].
    \end{aligned}
\end{equation}
We impose conditions on $c,\,k,\,\lambda\in(0,\bar\lambda)$ with $\bar\lambda$ as in Lemma~\ref{Lm:Exterior}, such that $u\leq \psi$ in
\begin{align}\label{eq_M}
    \mathcal{M}_{k,\lambda}=\left\{x\in \Omega : |x-y_\lambda|-d_\lambda\leq \frac{1}{k}\right\}\times I_\nu(t_o)\subset \Omega_{1/k}\times I_\nu(t_o).
\end{align}
\begin{center}
\begin{figure}
\psfragscanon
\psfrag{x}{$\scriptstyle x_o$}
\psfrag{y}{$\scriptstyle y$}
\psfrag{ylam}{$\scriptstyle y_\lm$}
\psfrag{Omega}{$\scriptstyle \Om$}
\psfrag{lam-lame}{$\scriptstyle \text{radius}=d_\lambda+\frac1k$}
\psfrag{ball}{$\scriptstyle B_{d_\lambda+\frac1k}(y_\lambda)$}
\psfrag{M}{$\scriptstyle {\mathcal M}_{k,\lambda}$}
\begin{center}
\includegraphics[width=.65\textwidth]{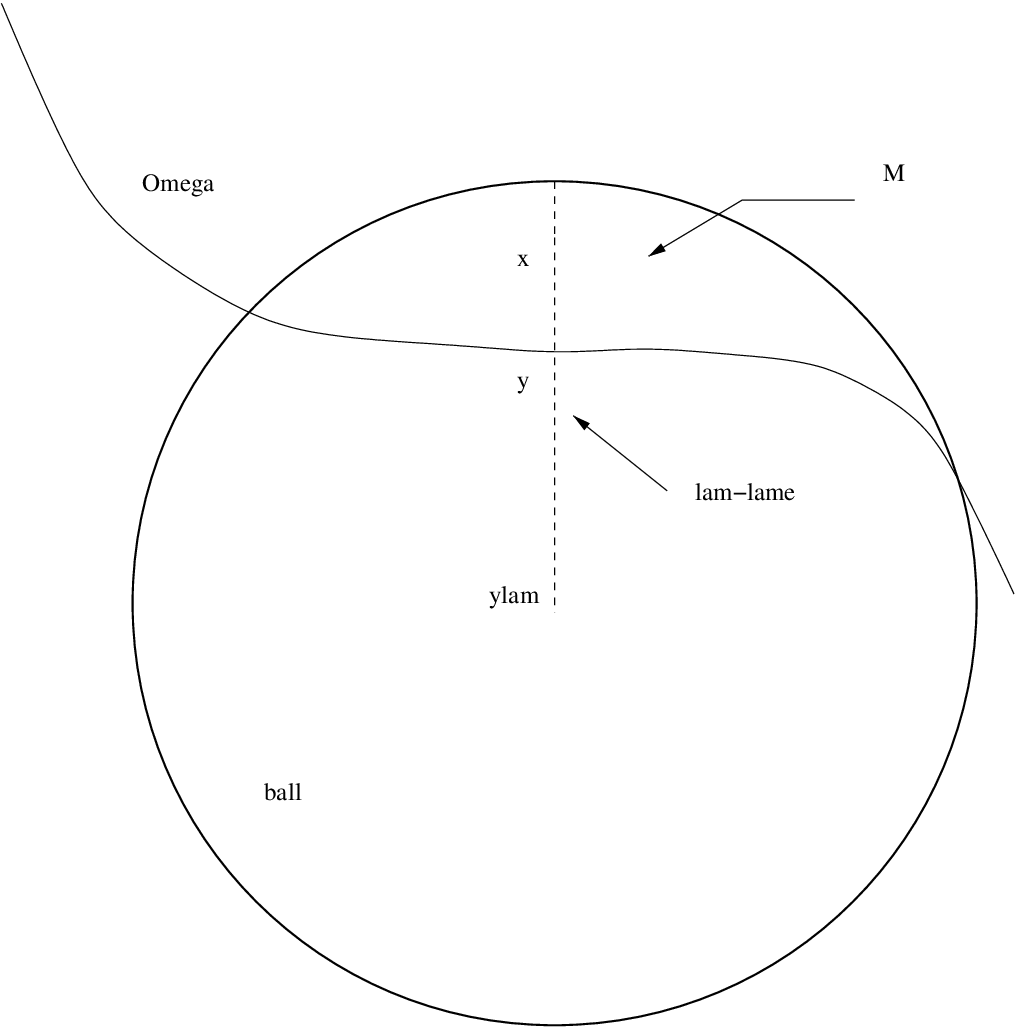}
\caption{\label{Fig:2}
}
\end{center}
\end{figure}
\end{center}
The set ${\mathcal M}_{k,\lambda}$ is described in Figure~\ref{Fig:2}. Since in the following we need $(x_o,t_o)\in{\mathcal M}_{k,\lambda}$, it is straightforward to check that by Lemma~\ref{Lm:Exterior} we must have $\frac1k>r$. In view of Theorem~\ref{Thm:Comp:Pr} it suffices to prove that
\[
\inf_{\partial_L\{\mathcal M_{k,\lm}\}}\psi>0,\quad \psi\Big|_{\partial_p\{\mathcal M_{k,\lm}\}}\ge u\Big|_{\partial_p\{\mathcal M_{k,\lm}\}},\quad\text{ and }\quad \partial_t\psi^q-\Delta_p \psi\ge0\,\,\text{ in }\,\,\mathcal{M}_{k,\lambda},
\]
where $\partial_L\{\mathcal M_{k,\lm}\}$ and $\partial_p\{\mathcal M_{k,\lm}\}$ denote respectively the lateral and parabolic boundaries of ${\mathcal M}_{k,\lm}$.

Note that clearly $\psi\geq 0$ in $\mathcal{M}_{k,\lm}$ and $\psi>0$ in $\partial_L\{\mathcal M_{k,\lm}\}$. We make the choice
\[
c=\max\left\{\left(1-e^{-1}\right)^{-1},\left(1-e^{-\nu^{\bar a+1}}\right)^{-1}\right\},
\]
so that in the set $\displaystyle \left\{x\in \Omega : |x-y_\lm|-d_\lm=\frac1k\right\}\times I_\nu(t_o)$
we have
\[
\psi\geq c M_\frac1k(1-e^{-1})\geq u.
\]
Moreover, in the set $\displaystyle \Omega_\frac{1}{k}\times\{t_o-\nu M^{q-p+1}\}$
it holds
\[
\eta\leq \exp\left(-\nu^{\bar a+1}\right)\implies \psi\geq cM_\frac{1}{k}\left(1-\exp\left(-\nu^{\bar a+1}\right)\right)\geq u.
\]
It remains to prove thar $\psi$ is a super-solution. First, note that
\[
\Delta_p \psi=(cM_\frac{1}{k}k)^{p-1}\eta^{p-1}\left(\frac{N-1}{|x-y_\lm|}-k(p-1)\right)
\]
and
\[
\frac{\partial}{\partial t}\psi^q=-q(\bar a +1)(cM_\frac{1}{k})^qM^{-(q-p+1)}(1-\eta)^{q-1}\eta \left(\frac{t_o-t}{M^{q-p+1}}\right)^{\bar a}.
\]
Therefore, we want to prove that in $\mathcal M_{k,\lm}$ 
\begin{align*}
    &(p-1)k^p-\frac{N-1}{|x-y_\lm|}k^{p-1}\\
    \geq\,&(\bar a+1)q(cM_\frac{1}{k})^{q-p+1}M^{-(q-p+1)}(1-\eta)^{q-1}\eta^{2-p}\left(\frac{t_o-t}{M^{q-p+1}}\right)^{\bar a}.
\end{align*}
Recall that 
\[
q-p+1>0,\quad M_\frac{1}{k}\leq M, \quad\inf_{x\in \Omega}|x-y_\lm|=d_\lm,\quad 0<\eta<1;
\]
hence, it suffices to prove
\begin{equation*}
    (p-1)k^p-\frac{N-1}{d_\lambda}k^{p-1}
    \geq\,(\bar a+1)qc^{q-p+1}(1-\eta)^{q-1}\eta^{2-p}\left(\frac{t_o-t}{M^{q-p+1}}\right)^{\bar a}.
\end{equation*}
We start by proving that $(1-\eta)^{q-1}\eta^{2-p}\left(\frac{t_o-t}{M^{q-p+1}}\right)^{\bar a}$ is bounded in $\mathcal M_{k,\lm}$. Since $\bar a=(1/q-1)_+$, we have to consider two alternatives.

\noindent\textit{\underline{I Alternative}.}\,\, If $q>1$ then $(1-\eta)^{q-1}\leq 1$ and in $\mathcal{M}_{k,\lm}$,
\[
\eta\geq e^{-1}\exp(-\nu^{\bar a+1}).
\]
If $p<2$ then $\eta^{2-p}\leq 1$; otherwise
\[
\eta^{2-p}\leq e^{-1}\exp(-\nu^{\bar a+1})^{2-p}.
\]
Setting
\[
C_1(\nu,p,q)\df=(p-1)^{-1}(\bar a+1)q\,c^{q-p+1}\max\left\{1,e^{-1}\exp(-\nu^{\bar a+1})^{2-p}\right\},
\]
we have to choose $k$, depending on $\lm$, such that
\[
k^p-\frac{N-1}{(p-1)d_\lm}k^{p-1}\geq C_1(\nu,p,q)
\]
which is definitely satisfied for
\[
k\geq \frac{N-1}{(p-1)d_\lm}+C_1(\nu,p,q)+1.
\]

\noindent\textit{\underline{II Alternative}.}\,\, If $0<q<1$ then $\bar a=1/q-1>0$ and
\begin{align*}
    \eta(x,t)\leq& \exp\left(-\left(\frac{t_o-t}{M^{p-q+1}}\right)^{\bar a+1}\right)\\
    \leq &1-\left(\frac{t_o-t}{M^{p-q+1}}\right)^{\bar a+1}+\left(\frac{t_o-t}{M^{p-q+1}}\right)^{2\bar a+2}.
\end{align*}
Hence, calling $\tau=\left(\frac{t_o-t}{M^{p-q+1}}\right)$ we get
\begin{align*}
    (1-\eta)^{q-1}\tau^{\bar a}\leq &\left(\tau^{1/q}-\tau ^{2/q}\right)^{q-1}\tau^{1/q-1}\\
    =& \left(1-\tau^{1/q}\right)^{q-1}\leq\left(1-\nu^{1/q}\right)^{q-1},
\end{align*}
since $\tau\in(0,\nu)$. Assuming $\nu<1$ we get $\left(1-\nu^{1/q}\right)^{q-1}\leq C_3(\nu,q)<\infty$. Noting that the bound on $\eta^{2-p}$ from the first alternative still holds, we have to choose $k$, depending on $\lm$, such that
\[
k^p-\frac{N-1}{(p-1)d_\lm}k^{p-1}\geq C_4(\nu,p,q)
\]
where $ C_4(\nu,p,q)=C_1(\nu,p,q)C_3(\nu,q) $. This is definitely satisfied for
\[
k\geq\frac{N-1}{(p-1)d_\lm}+C_4(\nu,p,q)+1.
\]
We satisfy the two alternatives by choosing 
\begin{align}\label{eq_kabove}
    k\geq\frac{N-1}{(p-1)d_\lm}+C(\nu,p,q)\quad\text{ for }\quad C=\max\{C_1,C_4\}+1.
\end{align}
Applying Theorem~\ref{Thm:Comp:Pr} and Lemma~\ref{Lm:Exterior} yields 
\begin{align*}
    u(x_o,t_o)\leq& \psi(x_o,t_o)\\
    \leq&cM_\frac{1}{k}k\,(d(x_o)+\lambda-d_\lambda)=c M_{\frac1k}\, k\, d(x_o).
\end{align*}
We still have to choose $k$, such that \eqref{eq_kabove} is satisfied, and we let  
\[
\frac{1}{k}=C\,r, \quad \lambda=d_\lambda=4C\frac{N-1}{p-1} r;
\]
hence, for $r<\bar r$ sufficiently small, depending on $C(\nu,p,q)$, \eqref{eq_kabove} holds. We have
\[
u(x_o,t_o)\le \frac cC M_{\frac1k}\frac{d(x_o)}{r}\le \gamma M \frac{d(x_o)}{r};
\]
since $(x_o,t_o)$ is arbitrarily chosen in $\Omega_r\times I_\nu$, possibly scaling back in order to take into account a general $\rho$, we conclude.
\hfill\os
\subsection{Control from below}

\begin{proposition}[Control from below]\label{Prop:below}
Let $0\le u\le M$ be a weak solution to \eqref{eq_doubly} in $\Om_T$, pick $s\in(0,T)$, and fix $x_o\in\partial\Om$. Assume that $u$ vanishes on $(\partial\Om\cap B_{2\rho}(x_o))\times(s,T)$.
For every $\nu>0$ there exist $\gamma_2>0$ and $\bar s\in(0,\frac12)$, both depending on $N, q, p, \nu$, and $\|\partial \Om\|_{1,1}$, such that for all $s+\nu M^{q-p+1}\rho^p<t<T$, and for all $x\in\Om\cap B_{2\bar s\rho(x_o)}$ with $\dist(x,\partial\Om)<\bar s\rho$ it holds
    \begin{align}\label{eq_lwbd}
        u(x,t)\geq \gamma_2 \left(\frac{\dist(x,\partial\Om)}\rho\right) \inf_{(\Om^{\bar s\rho}\cap B_{2\rho}(x_o))\times(s,T)} u.
    \end{align}
\end{proposition}
\noindent As already done in \S~\ref{SS:ctrl-above}, in the following we assume $\rho=1$. 
For $x\in\Om$, we let $d(x)$ and $\Gamma_x$ as before, and $\Omega^s:=\{x\in \Omega : d(x)>s\}$. Moreover, for $x\in\Om$ and $y\in \Gamma_x$, let
\begin{align*}
    y^\lambda\df=x-\lambda\frac{y-x}{|y-x|}, \quad d^\lambda=d(y^\lambda), \quad \lambda\geq 0.
\end{align*}
We have another simple geometrical result.
\begin{lemma}\label{Lm:Interior}
    There exist $\overline{s}, \overline{\lambda}>0$ depending only upon $\Om$, such that for every $x\in \Omega^{\bar{s}}\backslash\Omega^{2\bar{s}}$, $y\in \Gamma_x$ and $\lambda\in (0,\overline{\lambda})$ it holds
    \begin{align*}
        d(x)+\lambda=d^\lambda
    \end{align*}
\end{lemma}
\begin{center}
\begin{figure}
\psfragscanon
\psfrag{x}{$\scriptstyle x$}
\psfrag{y}{$\scriptstyle y$}
\psfrag{ylam}{$\scriptstyle y^\lm$}
\psfrag{Omega}{$\scriptstyle \Om$}
\psfrag{lam-lame}{$\scriptstyle \lambda\le\overline{\lambda}$}
\psfrag{ball}{$\scriptstyle B_{d^\lambda}(y^\lambda)$}
\psfrag{dx}{$\scriptstyle d(x)$}
\psfrag{dlam}{$\scriptstyle d^\lm=d(x)+\lm$}
\begin{center}
\includegraphics[width=.65\textwidth]{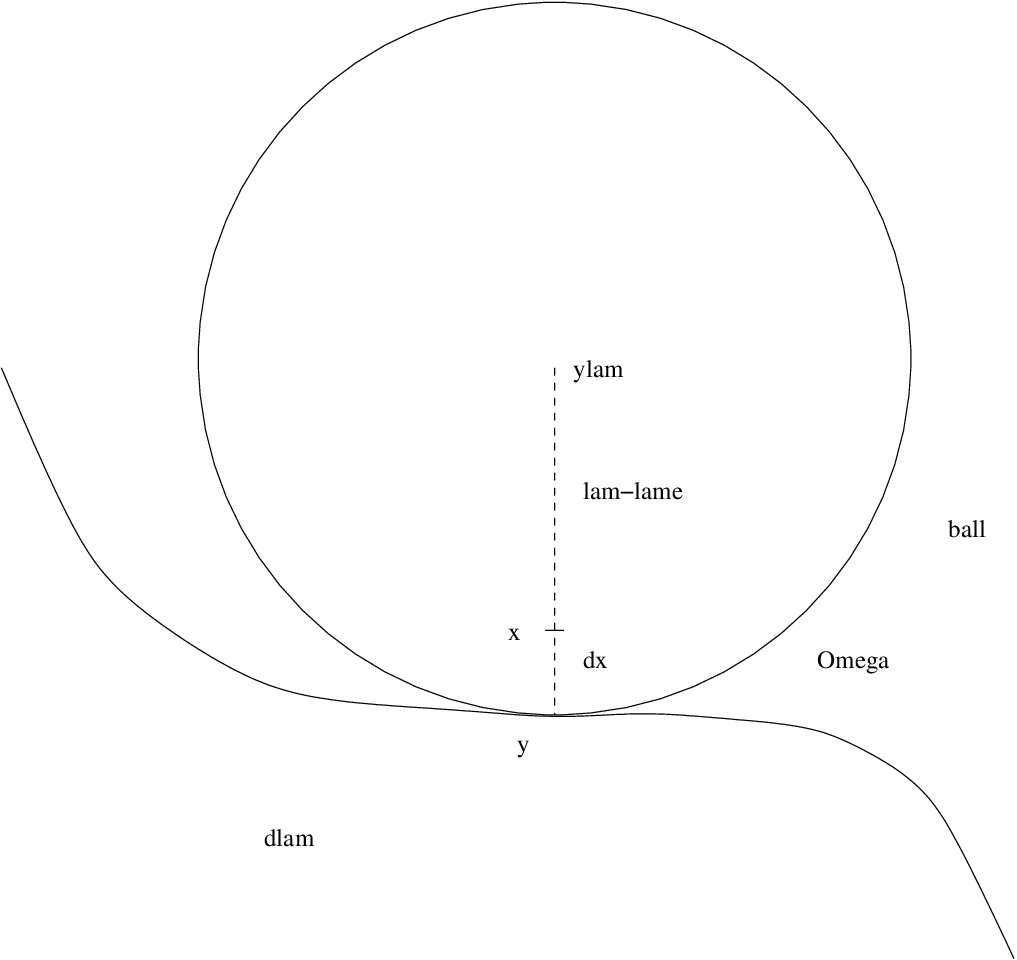}
\caption{\label{Fig:3}
}
\end{center}
\end{figure}
\end{center}
\begin{proof}
We take $\bar s=\frac{r_i}4$, as to ensure that for any $x\in\Om^{\bar s}\backslash\Om^{2\bar s}$ the set $\Gamma_x$ contains a single point due to Theorem~\ref{Thm:proj}. Moreover, by Definition~\ref{Def:UIBC}, it suffices to take $\overline{\lambda}=\frac{r_i}2$; indeed, for any $\lambda\in(0,\overline{\lambda})$, the corresponding $y^{\lambda}$ is the center of the ball of radius $d^\lambda=d(x)+\lambda<r_i$, that touches the boundary from the interior at $\Gamma_x$; by assumption, such a radius is uniform, and thus independent of $y$, and therefore also of $x$. Figure~\ref{Fig:3} describes what happens.
\end{proof}

In order to apply the Comparison Principle of Theorem \ref{Thm:Comp:Pr} with $u$ playing the role of super-solution, we have to prove that $u>0$ in the interior of $\Omega_T$. For this purpose, we start by recalling the following expansion of positivity (see \cite[Proposition 1.8]{BDGLS}).

\begin{proposition}\label{Prop:exp_pos}
    Let $u$ be a non-negative, local, weak super-solution to \eqref{eq_doubly} in $\Omega_T$. Suppose that for some constants $M>0$ and $\alpha\in (0,1)$, there holds
    \begin{align*}
        \big|\{u(\cdot,t_o)\geq M\}\cap K_\varrho(x_o)\big|\geq \alpha|K_\varrho|.
    \end{align*}
    Then, there exist constants $\delta,\, \eta\in (0,1)$, depending only on the data and $\alpha$, such that
    \[
    u\geq \eta M, \mbox{ a.e. in } K_{2\varrho}(x_o)\times \left(t_o+\frac{1}{2} \delta M^{q+1-p}\varrho^p,t_o+\delta M^{q+1-p}\varrho^p\right],
    \]
    provided
    \[
    K_{8\varrho}(x_o)\times \left(t_o,t_o+\delta M^{q+1-p}\varrho^p\right]\subset \Omega_T.
    \]
\end{proposition}
The previous statement is given in terms of cubes $K_\rho(x_o)$, but it holds for balls $B_\rho(x_o)$ as well. 
The proof of the next result follows the approach developed in \cite[Lemma~4.2]{moring2023supercaloric}.
\begin{lemma}\label{Lm:positive}
    Let $u$ be a continuous, non-negative, bounded, local, weak super-solution to \eqref{eq_doubly}, and assume $\Om$ bounded and connected. If
    \begin{align*}
       u(x_o,t)>0 \quad \mbox{ for some } x_o\in \Omega,\,\,t\in (0,T),
    \end{align*}
    then 
    \begin{align*}
       u(x,t)>0 \quad \mbox{ for any } x\in\Omega.
    \end{align*}
\end{lemma}
\begin{proof}
    Call $u_o:=u(x_o,t)>0$. We \textit{claim} that
    \[
    u(\cdot, t)>0 \quad \mbox{ in } B_r(x_o),
    \]
    for any $r$ such that $B_{16r}(x_o)\subset \Omega$. By continuity of $u$, we have that
\begin{align}\label{eq:uuo}
    u\geq \frac12 u_o \quad \mbox{ in } B_\rho(x_o)\times [t-\rho^p, t]\subset \Omega_T,
\end{align}
for some $\rho>0$. If $\rho>r$, we are done. Otherwise, apply 
Proposition~\ref{Prop:exp_pos} with $\alpha=(\rho/r)^N\in (0,1)$. Let $\delta=\delta(p,q,N,\alpha)$ be the number given by this proposition. Choose $M\in (0,u_o/2)$ so small that $\dsty\delta M^{q+1-p}r^p\leq \rho^p$, 
and let $\dsty t_o=t-\delta M^{q+1-p}r^p\in [t-\rho^p, t)$.
Because of \eqref{eq:uuo}, since $M\leq u_o/2$, and $\rho<r$, we have that
\[
|\{u(\cdot,t_o)\geq M\}\cap B_r(x_o)|\geq |B_\rho(x_o)|=\alpha |B_r(x_o)|.
\]
By the previous proposition,
\[
u(\cdot, t)\geq \eta M \quad \mbox{ in } B_{2r}(x_o),
\]
which proves our claim.

Next, we observe that by the previous proven claim if
\[
u(x_o,t)>0 \quad\mbox{ for some } x_o\in \Omega
\]
then
\[
u(x_1,t)>0 \quad\mbox{ for any } x_1\in \Omega.
\]
In order to show it, let $\Gamma\subset \Omega$ be a compact path connecting $x_o$ and $x_1$, which exists because $\Omega$ is bounded and connected. First of all, since $\Omega$ is open, there exists $\bar r>0$, depending only on $\Omega$, $x_o$, and $x_1$, and finitely many points $y_1,\, y_2,\, \hdots,\, y_n\in \Gamma$, such that $x_o=y_1,\, y_n=x_1$, $\Gamma$ is covered by the balls $B_{\bar r}(y_i)$, $B_{16\bar r}(y_i)\subset \Omega$ and $y_{i+1}\in B_{\bar r}(y_i)$. Indeed, if this were not the case, one could find sequences of points $(y_k)_k\subset \Gamma$ and radii $r_k\to 0$ such that $B_{r_k}(y_k)\not \subset \Omega$. Let $z_k\in B_{r_k}(y_k)\not \subset \Omega$. Since $\Gamma$ is compact, $y_k\to \bar y\in \Gamma$, up to a sub-sequence. Then
\begin{align*}
    |\bar y-z_k|\leq |\bar y-y_k|+|y_k-z_k|\leq \epsilon(k)+r_k\to 0\quad  \mbox{ as } k\to \infty,
\end{align*}
and therefore, $\bar y$ is not an interior point, which contradicts the fact that $\Omega$ is open.

To conclude the proof, we simply inductively apply the claim, starting with $y_1=x_o$ and get $u(\cdot,t)>0$ in $B_{\bar r}(y_1)$ which gives $u(y_2)>0$, and so on, eventually getting $u(y_n,t)>0$, as intended.
\end{proof}
\noindent\textit{Proof of Proposition~\ref{Prop:below}}. 
Let the interval $I_\nu$ be as before, and set
\begin{align*}
    \mu_r=\inf_{\Omega^r\times I_\nu} u(x,t).
\end{align*}

Choose $r\in(0,\frac{r_i}4)$. 
Let $(x_o,t_o)$ be a point in $(\Om^r\backslash\Om^{2r})\times I_\nu$, let $y\in\Ga_{x_o}$, and define
\[
I_\nu(t_o)\df=(t_o-\nu M^{q+1-p},t_o].
\]
We let $k>0$ and for $\bar\lambda$ as in Lemma~\ref{Lm:Interior}
\begin{align*}
&y^\lm=x_o-\lm\frac{y-x_o}{|y-x_o|},\quad d^\lm=d(y^\lm),\quad \lm\in(0,\overline{\lambda}),\\
&{\mathcal M}_{k,\lm}=\Big\{x\in\Om:\,d^\lm-\frac1k\le|x-y^\lm|\le d^\lm-\frac1{2k}\Big\}\times I_\nu(t_o).
\end{align*}
\begin{center}
\begin{figure}
\psfragscanon
\psfrag{x}{$\scriptstyle x_o$}
\psfrag{y}{$\scriptstyle y$}
\psfrag{ylam}{$\scriptstyle y^\lm$}
\psfrag{dx}{$\scriptstyle d(x)$}
\psfrag{Omega}{$\scriptstyle \Om$}
\psfrag{lam-lame}{$\scriptstyle \lm\le\overline{\lm}$}
\psfrag{ball}{$\scriptstyle B_{d^\lambda-\frac1{2k}}(y^\lambda)$}
\psfrag{ball2}{$\scriptstyle B_{d^\lambda-\frac1k}(y^\lambda)$}
\psfrag{M}{$\scriptstyle {\mathcal M}_{k,\lambda}$}
\begin{center}
\includegraphics[width=.65\textwidth]{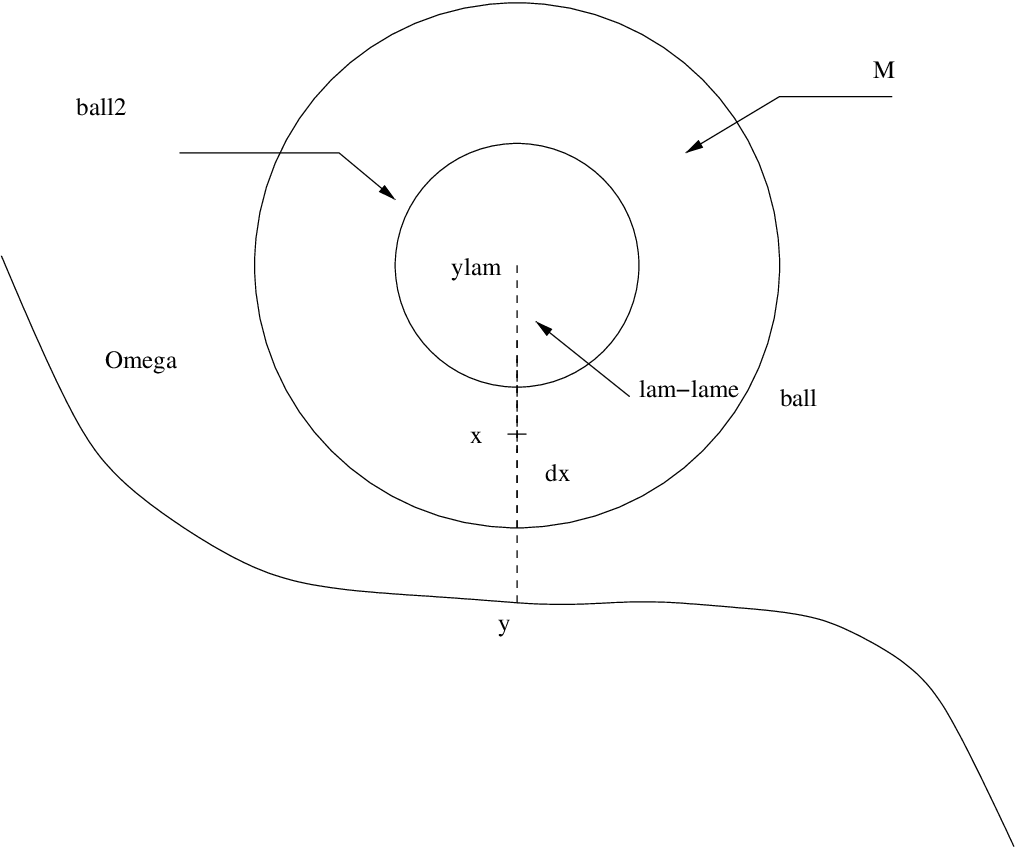}
\caption{\label{Fig:4}
}
\end{center}
\end{figure}
\end{center}
Since in the following we need $(x_o,t_o)\in{\mathcal M}_{k,\lm}$, it is straightforward to check that by Lemma~\ref{Lm:Interior} we must have $\frac1{2k}\le d(x_o)\le \frac1k$.
Let 
\[
\eta:{\mathcal M}_{k,\lm}\to\Rr,\quad\psi:{\mathcal M}_{k,\lm}\to\Rr
\]
be defined by
\begin{align*}
&\eta(x,t)=\exp\left[-k\left(|x-y^\lm|-d^\lm\right)\right]\exp\left[-\left(\frac{t_o-t}{\nu M^{q+1-p}}\right)^{\bar a+1}\right]\\
&\psi=\frac1e\mu_{\frac1k}(\eta-1)_+
\end{align*}
where $\bar a=(q^{-1}-1)^+$ as before. 

Note that by Lemma \ref{Lm:positive}, we have
\[
\essinf_{\partial_L \{{\mathcal M}_{k,\lm}\}}u>0,
\]
where we recall that $\partial_L$ denotes the lateral boundary of the cylindrical set. In the following, we want to choose $k$ and $\lm$ in such a way that 
\[
\psi\Big|_{\partial_p\{\mathcal M_{k,\lm}\}}\le u\Big|_{\partial_p\{\mathcal M_{k,\lm}\}},\quad\text{ and }\quad \partial_t\psi^q-\Delta_p \psi\le0\,\,\text{in}\,\,{\mathcal M}_{k,\lm},
\]
where $\partial_p\{\mathcal M_{k,\lm}\}$ denotes the parabolic boundary of ${\mathcal M}_{k,\lm}$. Since $u$ is a non-negative solution, by Theorem~\ref{Thm:Comp:Pr} we can then conclude that $\dsty\psi\le u$ a.e. in ${\mathcal M}_{k,\lm}$.

\noindent On the set $\dsty\Big\{|x-y^\lm|=d^\lm-\frac1{2k}\Big\}\times I_{\nu}(t_o)$ we have
\[
\eta(x,t)=e^\frac{1}{2}\cdot\exp\left[-\left(\frac{t_o-t}{\nu M^{q+1-p}}\right)^{\bar a+1}\right]\le e^\frac{1}{2}\,\,\Rightarrow\,\,\psi(x,t)\le\frac{e^\frac{1}{2}-1}{e}\mu_{\frac1k}\le u(x,t).
\]
On the set $\dsty\Big\{|x-y^\lm|=d^\lm-\frac1k\Big\}\times I_{\nu}(t_o)$ we have
\[
\eta(x,t)=e\cdot\exp\left[-\left(\frac{t_o-t}{\nu M^{q+1-p}}\right)^{\bar a+1}\right]\le e\,\,\Rightarrow\,\,\psi(x,t)\le\frac{e-1}{e}\mu_{\frac1k}\le u(x,t).
\]
Finally, at the bottom of the set, i.e. when $t=t_o-\nu M^{q+1-p}$, we have
\[
\eta(x,t)=e^{-1}\cdot\exp\left[-k\left(|x-y^\lm|-d^\lm\right)\right],
\]
so that $e^{-1}\le\eta\le e^{-\frac12}$, and $\psi=0$.
Now we determine what conditions on $k$ and $\lm$ yield that
$\psi$ is a sub-solution.
By straightforward computations, we have
\begin{align*}
&\partial_t \psi^q=\frac1{e^q}\,\mu_{\frac1k}^q q(\eta-1)_+^{q-1}\eta(\bar a+1)\left(\frac{t_o-t}{\nu M^{q+1-p}}\right)^{\bar a}\frac1{\nu M^{q+1-p}},\\
&\Delta_p\psi=-\left(\frac ke\,\mu_{\frac1k}\right)^{p-1}\chi_{[\eta>1]}\,\eta^{p-1}\left[\frac{N-1}{|x-y^\lm|}-k(p-1)\right].
\end{align*}
Hence, we want to know under what conditions
\begin{align*}
&\left(\frac ke\,\mu_{\frac1k}\right)^{p-1}\chi_{[\eta>1]}\,\eta^{p-1}\left[\frac{N-1}{|x-y^\lm|}-k(p-1)\right]\\
&+\frac1{e^q}\,\mu_{\frac1k}^q q(\eta-1)_+^{q-1}\eta(\bar a+1)\left(\frac{t_o-t}{\nu M^{q+1-p}}\right)^{\bar a}\frac1{\nu M^{q+1-p}}\le0.
\end{align*}
We can directly assume that $\eta\ge1$, otherwise there is nothing to prove. The previous inequality can be rewritten as
\begin{align*}
k^{p-1}\,\eta^{p-1}\left[k(p-1)-\frac{N-1}{|x-y^\lm|}\right]\ge\frac1{e^{q+1-p}}\,\frac{\mu_{\frac1k}^{q+1-p}}{M^{q+1-p}}\, \frac q\nu\,
(\eta-1)^{q-1}\eta^{2-p}(\bar a+1)\left(\frac{t_o-t}{\nu M^{q+1-p}}\right)^{\bar a}.
\end{align*}
Since $\bar a=(q^{-1}-1)^+$, we have to consider two alternatives.

\noindent\textit{\underline{I Alternative}.}\,\, If $q>1$, then $(\eta-1)^{q-1}\le\eta^{q-1}$, so that the right-hand side is bounded above by
\[
\frac1{e^{q+1-p}}\,\frac{\mu_{\frac1k}^{q+1-p}}{M^{q+1-p}}\, \frac q\nu\,
\eta^{q+1-p}\le\frac q\nu,
\]
and we end up estimating
\[
k^{p-1}\,\eta^{p-1}\left[k(p-1)-\frac{N-1}{|x-y^\lm|}\right]\ge\frac q\nu.
\]
\vskip.2truecm
\noindent\textit{\underline{II Alternative}.}\,\, If $q\le1$, then
\begin{align*}
&\eta(x,t)=\exp\left[-k\left(|x-y^\lm|-d^\lm\right)\right]\exp\left[-\left(\frac{t_o-t}{\nu M^{q+1-p}}\right)^{\frac1q}\right]\\
&\eta(x,t)-1\ge1\cdot\left[1+\left(\frac{t_o-t}{\nu M^{q+1-p}}\right)^{\frac1q}\right]-1
=\left(\frac{t_o-t}{\nu M^{q+1-p}}\right)^{\frac1q}.
\end{align*}
Hence,
\begin{align*}
(\eta-1)^{q-1}\le\left(\frac{t_o-t}{\nu M^{q+1-p}}\right)^{\frac{q-1}q}
\quad\Rightarrow\quad (\eta-1)^{q-1}(\bar a+1)\left(\frac{t_o-t}{\nu M^{q+1-p}}\right)^{\bar a}\le\frac1q.
\end{align*}
On the other hand, as we have already seen,
\[
e^{-1}\le\eta\le1\quad\Rightarrow\quad\eta^{2-p}\le C_1(p),\quad \eta^{p-1}\ge C_2(p),
\]
independently whether $p>2$ or $1<p<2$. We conclude that in any case
\[
k^{p-1}\left[k(p-1)-\frac{N-1}{|x-y^\lm|}\right]\ge\frac {C_3(p)}\nu.
\]
Provided we choose $k>1$ sufficiently large, which we can always do, it suffices to have
\[
k(p-1)-\frac{N-1}{|x-y^\lm|}\ge C_4(p,q,\nu),
\]
and also 
\begin{equation}\label{Eq:k-lb2-1}
k\ge C_5(p,q,\nu)+\frac{N-1}{(p-1)|x-y^\lm|}.
\end{equation}
Since we are working in the set $\dsty \Big\{d^\lm-\frac1k\le|x-y^\lm|\le d^\lm-\frac1{2k}\Big\}$, it is apparent that
\[
\frac{1}{|x-y^\lm|}>\frac1{d^\lm};
\]
hence,
\[
\partial_t\psi^q-\Delta_p\psi\le0\quad\text{ in }\,\,{\mathcal M}_{k,\lm},
\]
provided 
\begin{equation}\label{Eq:k-lb2-2}
\dsty k\ge C_5+\frac{(N-1)}{(p-1)d^\lm}.
\end{equation}
Theorem~\ref{Thm:Comp:Pr} yields $\dsty u\ge\psi$ in $\dsty{\mathcal M}_{k,\lm}$,
and in particular, $\dsty u(x_o,t_o)\ge\psi(x_o,t_o)$.
By the sheer definition of $\psi$, it is straightforward to conclude that
\[
u(x_o,t_o)\ge\frac1e\, \mu_{\frac1k}\, k[d^\lm-\lm]=\frac1e\, \mu_{\frac1k}\,k\,d(x_o).
\]
We can now choose $k$ and $\lm$ in terms of $r$, provided \eqref{Eq:k-lb2-2} is satisfied. We let
\[
d^\lm=\frac{3(N-1)}{p-1}r,\quad \frac1k=2r, 
\]
and so for $r< r_i$ sufficiently small, depending on $\nu,\,p,\,q$, \eqref{Eq:k-lb2-2} holds. We have
\[
u(x_o,t_o)\ge \gamma\frac{d(x_o)}{r}\mu_{2r},
\]
and since $(x_o,t_o)$ is arbitrarily chosen in $(\Om^r\backslash\Om^{2r})\times I_\nu$, we conclude. \hfill $\square$

\section{Boundary Harnack Inequality}\label{S:BHI}
For the next result we use some notation introduced in Section \ref{Sec:Carl}.

\begin{theorem}\label{Thm:BdHarn}
    Let $\Om$ be a $C^{1,1}$ domain, $u,\,v$ two weak solutions to \eqref{eq_doubly} in $\Om_T$ which satisfy the assumptions of Theorem~\ref{Thm:Carleson}, and such that $0<u,\,v\leq M$. Then, there exist positive constants $\bar s$ and $\gamma$, depending on $N$, $p$, $q$, $\|\partial \Om\|_{1,1}$, $\rho_o$, 
    such that
    \begin{equation}\label{eq_BdHarn}
            \gamma^{-1}\frac{\inf_{\tau \in I(t_o,\rho,2M^{q+1-p})}u(x_\rho,\tau)}{\sup_{\tau \in I(t_o,\rho,2M^{q+1-p})}v(x_\rho,\tau)}\leq \frac{u(x,t)}{v(x,t)}
        \leq \gamma\,\frac{\sup_{\tau \in I(t_o,\rho,2M^{q+1-p})}u(x_\rho,\tau)}{\inf_{\tau \in I(t_o,\rho,2M^{q+1-p})}v(x_\rho,\tau)}
    \end{equation}
    for all $(x,t)\in \{x\in K_{\bar s \rho/4}(x_o)\cap \Omega \,:\, \dist(x,\partial \Omega)<\bar s\rho /8\}$ with $\rho<\rho_o$ and $(x_o,t_o)\in S_T$.
    
\end{theorem}

\begin{proof}
    The proof follows an argument similar to the one developed in \cite[Theorem 3.8]{AGS}.

    Let $\rho\in (0,r_o)$ where $r_o$ was introduced in Definition \ref{Def:Lip:Dom}. In Propositions \ref{Prop:above} and \ref{Prop:below}, let
\[
    s=t_o-M^{q+1-p}\rho^p,\quad T=t_o+M^{q+1-p}\rho^p,
\]
and $\nu=1/2$. Define
\[
    V_{\frac{1}{2}, \rho}(x_o,t_o)=\{x\in K_{\bar s \rho /4}(x_o)\, : \, \dist(x,\partial\Omega)<\bar \rho/8\}\times I(t_o,\rho,\frac{1}{2}M^{q+1-p}),
\]
where $\bar s$ is given by Propositions~\ref{Prop:above} and \ref{Prop:below} when $\nu=1/2$. Applying these propositions to both $u$ and $v$ we can write, for every $(x,t)\in V_{1/2,\rho}(x_o,t_o)$,
\begin{equation}\label{eq_6.1}
    \begin{aligned}
        & \gamma_2 \dist(x,\partial \Omega)\,\mu_u(\bar s\rho)\leq \rho\, u(x,t)\leq \gamma_1 \dist(x,\partial \Omega)\,M_u(\rho),\\
        &\gamma_2 \dist(x,\partial \Omega)\,\mu_v(\bar s\rho)\leq \rho\, v(x,t)\leq \gamma_1 \dist(x,\partial \Omega)\,M_v(\rho),
    \end{aligned}
\end{equation}
with
\begin{align*}
    M_u(\rho)=\sup_{(\Omega\cap K_{\rho/4}(x_o))\times I(t_o,\rho,M^{q+1-p})}u,\quad M_v(\rho)=\sup_{(\Omega\cap K_{\rho/4}(x_o))\times I(t_o,\rho,M^{q+1-p})}v,
\end{align*}
and
\begin{align*}
    \mu_u(\bar s\rho)=\sup_{(\Omega^{\bar s\rho/8}\cap K_{\rho/4}(x_o))\times I(t_o,\rho,M^{q+1-p})}u, \quad \mu_v(\bar s\rho)=\sup_{(\Omega^{\bar s\rho/8}\cap K_{\rho/4}(x_o))\times I(t_o,\rho,M^{q+1-p})}v,
\end{align*}
where $\Omega^s:=\{x\in \Omega \, : \, s/2\leq \dist(x,\partial \Omega)\leq 2s\}$. 
We recall that by 
Propositions~\ref{Prop:above} and \ref{Prop:below}, $\gamma_1$ and $\gamma_2$ depend also on $\|\partial\Om\|_{C^{1,1}}$. Then, \eqref{eq_6.1} gives 
\begin{align*}
    \frac{\gamma_2}{\gamma_1}\frac{\mu_u(\bar s\rho)}{M_v(\rho)}\leq \frac{u(x,t)}{v(x,t)}\leq \frac{\gamma_1}{\gamma_2}\frac{M_u(\rho)}{\mu_v(\bar s\rho)}.
\end{align*}
By the Carleson estimate in Theorem~\ref{Thm:Carleson}, we can bound above $M_u(\rho)$ by
\[
    M_u(\rho)\leq \gamma_1\sup_{\tau\in I(t_o,\rho, 2M^{q+1-p})}u(x_\rho,\tau).
\]
On the other hand, $\mu_u(\bar s\rho)$ is attained at some point $(x_*,t_*)\in K_{\rho/4}(x_o)\times I(t_o,\rho, M^{q+1-p})$. By the Harnack inequality of Theorem \ref{Thm:Harnack},
\[
    \gamma_3^{-1}u(x_\rho,t_*)\leq u(x_*,t_*)\leq \gamma_3 u(x_\rho,t_*)
\]
where $\gamma_3$ depends on $N$, $p$, $q$ (but it depends neither on $\|\partial \Om\|_{1,1}$, nor on $\bar s$), and therefore
\[
\mu_u(\bar s \rho)\geq \gamma_3^{-1}\inf_{\tau \in I(t_o,\rho, 2M^{q+1-p})} u(x_\rho,\tau).
\]
In an identical way we obtain analogous estimates concerning $v$. Combining everything, we get \eqref{eq_BdHarn} for all $(x,t)\in V_{\frac{1}{2}, \rho}(x_o,t_o)$.
\end{proof}

\appendix
\appendixpage
\section{Introduction}\label{S:A:Intro}
We are interested in the H\"older continuity of non-negative, local, bounded solutions to the doubly nonlinear equation \eqref{eq:structure}--\eqref{Def:structure}.
Moreover, we assume that $\partial\Omega$ satisfies 
the so-called property of {\it positive geometric density}, that is, there 
exist $\be\in(0,1)$ and $R_*>0$ such that for all 
$y\in\partial\Omega$, and for all  $0<\rho\le R_*$
\begin{equation}\label{Eq:11:8:4}
\big|B_{\rho}(y)\cap(\rn\backslash\Omega)\big|\ge \be\big|B_\rho\big|.
\end{equation} 
Without loss of generality, \eqref{Eq:11:8:4} can be stated in terms of cubes, instead of balls.
Finally, we suppose that $u$ vanishes at the lateral boundary $S_T$ in the sense of Definition~\ref{def-lateral-trace}. We are not interested either in non-homogeneous Dirichlet boundary conditions, nor on initial conditions.

It is apparent that \eqref{eq_doubly} is just the prototype case of the wider class of equations we study here; we can work in this more general framework with no major extra effort, since the methods we employ are purely structural, and are not limited to particular properties of the prototype case. Moreover, we think that our results may have a broader interest, beyond the specific problem we are dealing with here.

We proceed in the following way:
\begin{enumerate}
\item We prove the H\"older continuity of non-negative, bounded solutions in the interior of $\Omega_T$; more precisely, we show that the following Geometric Decay Lemma holds:
\begin{lemma}\label{Lm:Geom-Decay}
\emph{(H\"older continuity)} Let $u:\Om_T\to\Rr$ be a non-negative, locally bounded, local, weak solution to \eqref{eq:structure}--\eqref{Def:structure} and $(x_o,t_o)\in\Om_T$. There exist constants $\varepsilon$,
$\delta\in(0,1)$, and $A\ge1$, 
depending only on $p$, $q$, $N$, $C_o$, $C_1$ such that, setting
\begin{equation*}
\om_n = \delta\om_{n-1}, \quad
\theta_n = \left ( \frac{\om_n}A\right )^{q+1-p}, \quad
\rho_n = \varepsilon \rho_{n-1} \quad \textrm{and} \quad
Q_n = (x_o,t_o)+Q_{\rho_n}^- (\theta_n),
\end{equation*}
for all $n\in\Nn_0$, there holds $Q_{n+1}\subset Q_n\subset\dots\subset\subset\Om_T$ and
\begin{equation}\label{Eq:geom-decay}
\essosc_{Q_n} u \le \om_n.
\end{equation}
\end{lemma}

Once such a lemma is available, the local H\"older continuity of $u$ in the interior follows quite straightforwardly.
\item We prove the H\"older continuity of non-negative solutions up to the lateral boundary of $\Omega_T$, provided that $\partial\Omega$ satisfies a positive geometric density condition, and $u$ vanishes at the lateral boundary; in particular, if $\partial\Omega$ is flat at $x_o$, as explained in Section~\ref{Sec:Carl}, it satisfies condition \eqref{Eq:11:8:4} with $\beta=\frac12$ at such a point.
Once more, we do not limit ourselves to the result we need in this particular framework, but we deal with a broader class of boundaries, since it does not pose further difficulties. 

As in the previous step, the crucial point lies in proving a Geometric Decay Lemma of the solution $u$ at the boundary. With respect to Lemma~\ref{Lm:Geom-Decay}, there are only two (obvious) differences, namely $(x_o,t_o)\in S_T$, and \eqref{Eq:geom-decay} is substituted by
\begin{equation}\label{Eq:geom-decay-bd}
\essosc_{Q_n\cap\Om_T} u \le \om_n.
\end{equation}
As before, the H\"older continuity of $u$ is a straightforward consequence of \eqref{Eq:geom-decay-bd}.
\end{enumerate}
Once these first two tasks are accomplished, we rely on them to study the H\"older continuity of a particular class of weak, signed solutions to the doubly nonlinear equation
\begin{equation}\label{Eq:A1:1bis}
\partial_t\left(|u|^{q-1}u\right)-\dvg \mathbb{A}(x, t, u, D u)=0,
\end{equation}
where $\mathbb{A}$, $p$, $q$ satisfy the same assumptions as in \eqref{Def:structure}. We have two further steps:
\begin{enumerate}
\setcounter{enumi}{2}
\item We prove that by odd reflection along a flat boundary, non-negative solutions to \eqref{eq:structure} become signed solutions to \eqref{Eq:A1:1bis};
\item We prove that the signed solutions we have constructed in this manner continue to satisfy the Geometric Decay Lemma. This controlled decay is precisely what we need for the proof of the Carleson estimates of Theorem~\ref{Thm:Carleson}. 
\end{enumerate}
As a consequence of these four steps, we can conclude that the signed solutions constructed by odd reflection of non-negative solutions along a flat boundary satisfy
$$
\omega_{n+1}=\delta \omega_n\quad\text{ for a proper }\,\delta\in(0,1),\quad\text{ and }\quad \essosc_{Q_n} u\le\om_n,
$$
independently of the kind of cylinder we use to define $\omega_n$ and $\dsty\essosc_{Q_n} u$. We refrain from providing a complete and formal statement of the Geometric Decay Lemma in this final situation, since it is a trivial adaptation of Lemma~\ref{Lm:Geom-Decay}.

\section{H\"older Continuity in the Interior}
We assume $u$ to be a non-negative, local, bounded, weak solution in $\Omega_T$. We closely follow \cite{VV_2022}. 
\subsection{Preliminary material}

Here we include preliminary material which is used in the next sections. These results can all be found in \cite{BDGLS}.

The following functions $\mathfrak{g}_\pm$ are instrumental in our argument. For $w, k\in \Rr$ and $q>0$ we define 
\begin{align}\label{eq:g_function}
    \mathfrak{g}_\pm(w,k):=\pm q\int_k^w |s|^{q-1} (s-k)_\pm dx.
\end{align}
We have the following lemma.
\begin{lemma}\label{Lm:g}
    Let $q>0$ and $\mathfrak{g}_\pm$ be defined by \eqref{eq:g_function}. There exists a constant $\gamma(q)$ such that for all $a, b\in \Rr$,
    \[
        \frac{1}{\gamma}\left(|a|+|b|\right)^{q-1}(a-b)^2_\pm\leq \mathfrak{g}_\pm(a,b)\leq \gamma \left(|a|+|b|\right)^{q-1}(a-b)^2_\pm
    \]
\end{lemma}
We recall the energy estimate satisfied by $\mathfrak{g}_\pm$.
\begin{proposition}\label{Prop:energy}
    Let $p>1$ and $q>0$. There exists a constant $\gamma(C_o, C_1, p)>0$, such that whenever $u$ is a non-negative weak sub(super)-solution to \eqref{def-eq-global}, $Q_{\rho, s}= K_\rho(x_o)\times (t_o-s, t_o]\Subset \Omega_T$ is a parabolic cylinder, $k>0$, and $\eta$ any non-negative, piecewise smooth cutoff function vanishing on $\partial K_\rho(x_o)\times (t_o-s,t_o)$, we have
    \begin{align*}
        &\max\left\{ \esssup_{t_o-s<t<t_o}\int_{K_\rho(x_o)\times \{t\}} \eta^p \mathfrak{g}_\pm(u,k)\, dx\, , \, \iint_{Q_{\rho,s}}\eta^p|D(u-k)_\pm|^p dxdt   \right\}\\
        &\quad \leq \, \gamma \iint_{Q_{\rho,s}}\left[ (u-k)^p_\pm |D\eta|^p+\mathfrak{g}_\pm(u,k) |\partial_t\eta|  \right]\, dxdt\\
        &\qquad +\int_{K_\rho(x_o)\times\{t_o-s\}} \eta^p\mathfrak{g}_\pm(u,k)\,dx.
    \end{align*}
\end{proposition}
We conclude with the following geometric convergence result.
\begin{lemma}\label{Lm:fast_convergence}
    Let $(Y_n)_{n\in \Nn_o}$ be a sequence of positive real numbers satisfying the recursive inequalities
    \[
        Y_{n+1}\leq Cb^n Y_n^{1+\delta},
    \]
    where $C,b>1$ and $\delta>0$ are given numbers. If
    \[
        Y_o\leq C^{-1/\delta}b^{-1/\delta^2},
    \]
    then $Y_n\to 0$ as $n\to \infty$.
\end{lemma}
\subsection{A De Giorgi-type Lemma}\label{SS:A:B2}
For $\delta \in(0,1)$ to be fixed later, let $R_o>0$ be so small that
\begin{equation}\label{Eq:A1:3}
Q_o^*\df=K_{R_o}(x_o) \times(t_o-\delta R_o^p, t_o] \subset \Omega_T,
\end{equation}
and let
$$
\tilde{\mu}_o^+= \esssup _{Q_o^{*}} u,\quad\tilde{\mu}_o^-= \essinf_{Q_o^{*}} u, \quad\tilde{\omega}_o=\essosc_{Q_o^*} u=\tilde{\mu}_o^{+}-\tilde{\mu}_o^{-}.
$$
If $\tilde{\mu}_o^->0$, we are led to consider the H\"older continuity of a $p$-parabolic-like equation. 
We thoroughly discuss this possibility later on, at the end of the section. For the moment, let us assume that $\tilde{\mu}_o^{-}=0$. Hence,
$$
\tilde{\mu}_o^{+}=\tilde{\omega}_o.
$$
Since $u$ is locally bounded, without loss of generality we may assume that $\tilde{\omega}_o \le 1$, so that
$$
Q_o\df=K_{R_o}(x_o) \times(t_o- \tilde{\omega}_o^{q+1-p} \delta R_o^p, t_o] \subseteq Q_o^* \subset \Omega_T.
$$
We let
$$
\mu_o^{+}=\esssup_{Q_o} u,\quad\mu_o^{-}=\essinf_{Q_o} u.
$$
Since $\mu_o^{+} \le \tilde{\mu}_o^{+}$, it is apparent that if we let
$$
Q_1\df=K_{R_o}(x_o) \times(t_o-\left(\mu_o^{+}\right)^{q+1-p} \delta R_o^p, t_o],
$$
then $\dsty Q_1 \subseteq Q_o$. 
We have the following.
\begin{lemma}\label{Lm:A1:1} 
Let $u$ be a non-negative, local, bounded, weak sub-solution to \eqref{eq:structure}--\eqref{Def:structure} in $\Omega_T$ for 
$0<p-1<q<\frac{N(p-1)}{(N-p)_{+}}$. Set
$\dsty\theta=\delta(\mu_o^{+})^{q+1-p}$
for the same $\delta \in(0,1)$ as in \eqref{Eq:A1:3}. Assume that $\dsty (x_o,t_o)+Q_\rho(\theta) \subset Q_1$ and
\begin{equation}\label{Eq:A1:4}
\esssup_{(x_o,t_o)+Q_\rho(\theta)} u \le \mu_o^{+}.
\end{equation}
Then, there exists a constant $\nu>0$ depending only on the data $\left\{q, p, N, C_o, C_1\right\}$, and $\delta$, such that if
\begin{equation}\label{Eq:A1:5}
\left|\left\{u>\frac{\mu_o^{+}}{2}\right\} \cap(x_o,t_o)+Q_\rho(\theta)\right| \le \nu\left|Q_\rho(\theta)\right|
\end{equation}
then
\begin{equation}\label{Eq:A1:6}
u<\frac{3}{4} \mu_o^{+} \quad \text { a.e. in } \,(x_o,t_o)+Q_{\frac{\rho}{2}}(\theta).
\end{equation}
Moreover, for some $\nu_o$ depending only on the data $\left\{q, p, N, C_o, C_1\right\}$, 
we have $\dsty\nu=\nu_o \delta^{N/p}$.
\end{lemma}
\begin{proof}
Without loss of generality, we may assume $(x_o,t_o)\equiv(0,0)$. In this proof, $c$ always denotes a constant depending on $N, p, q, C_o, C_1$, which may change from line to line.

Define sequences of numbers and sets as follows
\begin{align*}
    \rho_j:=\frac{1}{2}\left(1+\frac{1}{2^j}\right),\, k_j:=\left(1-\frac{1}{4}-\frac{1}{2^{j+2}}\right)\mu_+,\\
    Q_j:=Q_{\rho_j,\theta \rho_j^p}, \, A_j:=Q_j\cap \{u>k_j\} ,\, Y_j:=\frac{|A_j|}{|Q_j|}.
\end{align*}
Let $\varphi\in C^\infty(Q_j;[0,1])$ be a function vanishing near the parabolic boundary of $Q_j$, satisfying $\varphi\equiv1$ on $Q_{j+1}$, and such that
\[
    |D\varphi_j|\leq \rho^{-1}2^{j+2},\, |\partial_t \varphi_j|\leq c\,\theta^{-1}\rho^{-p}2^{jp}.
\]
Note that in the set where $u>k_{j+1}$ we have 
\begin{align}\label{eq:B8.3}
    u-k_j>k_{j+1}-k_j=\frac{\mu_+}{2^{j+3}}.
\end{align}
Together with an application of H\"older inequality this yields
\begin{align}\label{eq:B8.4}
    \frac{\mu_+^p}{2^{(j+3)p}}|A_{j+1}|\leq&\, \iint_{A_j+1}(u-k_j)_+^p\,dxdt\\
    \leq &\,\left[\iint_{A_j+1}(u-k_j)_+^{p\frac{N+p}{p}}\,dxdt\right]^\frac{N}{N+p}|A_{j+1}|^\frac{p}{N+p}.\nonumber
\end{align}
Using again \eqref{eq:B8.3} and the H\"older inequality we further get
\begin{align*}
    &\iint_{A_j+1}(u-k_j)_+^{p\frac{N+p}{p}}\,dxdt\leq \int_{T_{j+1}}\int_{B_{j+1}}(u-k_j)_+^\frac{p^2}{N}\,dxdt\\
    \leq\,&\int_{T_{j+1}}\left[\int_{B_{j+1}}(u-k_j)_+^p\chi_{A_{j+1}}\,dx\right]^\frac{p}{N}\left[\int_{B_{j+1}}(u-k_j)_+^{p^*}\,dx \right]^\frac{p}{p^*}dt.
\end{align*}
We used $Q_j=B_j\times T_j$, where $B_j=B_{\rho_j}$ and $T_j=(-\theta\rho_j^p,0)$. Writing
\[
    (u-k_j)_+^p=(u-k_j)_+^{p-1-q}(u-k_j)_+^{q+1}.
\]
and since $p-1-q<0$, by \eqref{eq:B8.3} we estimate in $A_{j+1}$,
\[
    (u-k_j)_+^{p-1-q}<\left(\frac{\mu_+}{2^{j+3}}\right)^{p-1-q}.
\]
Consider first the case $q\geq1$: by Lemma \ref{Lm:g},
\begin{align*}
    (u-k_j)_+^{q+1}\leq &\, (u-k_j)_+^2(|u|+|k_j|)^{q-1}\\
    \leq &c\,  \mathfrak{g}_+(u,k_j).
\end{align*}
Hence in $A_{j+1}$ 
\begin{align*}
    (u-k_j)^p_+\leq c \left(\frac{\mu_+}{2^{j}}\right)^{p-1-q} \mathfrak{g}_+(u,k_j).
\end{align*}
Consider the case $0<q<1$: then,
\begin{align*}
    \mathfrak{g}_+(u,k_j)=&\,q\int_{k_j}^u (s-k_j)_+\,ds
    \geq q u^{q-1}\int_{k_j}^u(s-k_j)_+\,ds \\
    =&\, \frac{q}{2}u^{q-1}(u-k_j)^2_+    =\, \frac{q}{2}\left(\frac{(u-k_j)_+}{u}\right)^{1-q}(u-k_j)_+^{q+1}.
\end{align*}
Since $f(t)=t/(t-a)$ is monotone decreasing in $(a,+\infty)$, in the set $\{u>k_{j+1}\}$ it holds 
\[
\frac{u}{u-k_j}\leq \frac{k_{j+1}}{k_{j+1}-k_j}\leq 2^{j+3} .
\]
Thus,
\begin{align*}
    &(u-k_j)_+^{q+1}\leq \frac{2}{q}2^{(1-q)(j+3)}\mathfrak{g}_+(u,k_j)\leq  c\, 2^{j}\mathfrak{g}_+(u,k_j),
\end{align*}
and in $A_{j+1}$ we have
\[
    (u-k_j)_+^p\leq c \left(\frac{\mu_+}{2^{j}}\right)^{p-1-q} 2^j\mathfrak{g}_+(u,k_j).
\]
Therefore, we can combine both cases and get that for every $q>0$, in $A_{j+1}$ it holds
\[
    (u-k_j)_+^p\leq c \mu_+^{(p-1-q)} 2^{j(2+q-p)}\mathfrak{g}_+(u,k_j).
\]
Now we introduce the cutoff function $\varphi_j$, after which we can apply the energy estimate in Proposition \ref{Prop:energy} and Sobolev inequality. Thus,
\begin{align*}
    & \iint_{A_{j+1}} (u-k_j)_+^{p\frac{N+p}{N}}\,dxdt\\
    \leq &\, c \left(\frac{\mu_+^{(p-1-q)}}{2^{j(p-2-q)}}\right)^\frac{p}{N}\left[ \esssup_{T_j}\int_{B_j}\mathfrak{g}_+(u,k_j)\varphi^p_j\,dx \right]^\frac{p}{N}\int_{T_j}\left[\int_{B_j}[(u-k_j)_+\varphi_j]^{p^*}dx \right]^\frac{p}{p^*}dt\\
    \leq &\, c \left(\frac{\mu_+^{(p-1-q)}}{2^{j(p-2-q)}}\right)^\frac{p}{N}\left[ \esssup_{T_j}\int_{B_j}\mathfrak{g}_+(u,k_j)\varphi^p_j\,dx \right]^\frac{p}{N}\int_{T_j}\int_{B_j}|D[(u-k_j)_+\varphi_j]|^pdx dt\\
    \leq &\,c \left(\frac{\mu_+^{(p-1-q)}}{2^{j(p-2-q)}}\right)^\frac{p}{N}\left[ \iint_{A_j}\left[(u-k_j)_+^p|D\varphi_j|^p+\mathfrak{g}_+(u,k_j)|\partial_t\varphi| \right]dxdt \right]^\frac{N+p}{N}
\end{align*}
By Lemma \ref{Lm:g}
\begin{align*}
    \mathfrak{g}_+(u,k_j)\leq&\, c \left(|u|+|k_j|\right)^{q-1}(u-k_j)^2_+\\
    \leq &\, c 2^{q-1}|u|^{q-1}\leq c 2^{q-1}\mu_+^{q+1}.
\end{align*}
Recall also the estimates for the cutoff function
\[
|\varphi_j|\leq 1,\quad |D\varphi_j|\leq \rho^{-1}2^{j+2},\quad
|\partial_t \varphi_j|\leq c\theta^{-1}\rho^{-p}2^{jp}.
\]
Combining everything yields
\begin{align*}
    &\iint_{A_{j+1}} (u-k_j)_+^{p\frac{N+p}{N}}\,dxdt\\
    \leq &\,c \left(\frac{\mu_+^{(p-1-q)}}{2^{j(p-2-q)}}\right)^\frac{p}{N} \left[\left( \mu_+^p \rho^{-p} 2^{(j+2)p}+\varphi 2^{q-1} \mu_+^{q+1}c \theta^{-1} 2^{jp} \right)|A_j| \right]^\frac{N+p}{N}.
\end{align*}
Using \eqref{eq:B8.4} and the fact that $|A_{j+1}|\leq |A_j|$, we get
\begin{align*}
    |A_{j+1}|&\\
    \leq&\,c \left(\frac{\mu_+^{(p-1-q)}}{2^{j(p-2-q)}}\right)^\frac{p}{N+p} \left( \mu_+^p\rho^{-p}2^{(j+2)p}+c 2^{q-1} \mu_+^{q+1} \theta^{-1}2^{jp} \right) |A_j|^{1+\frac{p}{N+p}} 2^{(j+3)p}\mu_+^{-p}\\
    \leq &\,c \left(\frac{\mu_+^{(p-1-q)}}{2^{j(p-2-q)}}\right)^\frac{p}{N+p} \left( \rho^{-p}2^{2p}+c 2^{q-1} \mu_+^{q+1-p} \theta^{-1} \right) |A_j|^{1+\frac{p}{N+p}} 2^{(j+3)p}.
\end{align*}
Recall that $\theta^{-1}\mu_+^{q+1-p}=\varepsilon^{-1}$, $2p<2q+2$ and $\rho<1$. Thus
\begin{align*}
    |A_{j+1}|\leq &\,c \left(\frac{\mu_+^{(p-1-q)}}{2^{j(p-2-q)}}\right)^\frac{p}{N+p}2^{2jp+3p}\rho^{-p}2^{q-1}\left(2^{q+3}+c\varepsilon^{-1}\right)|A_j|^{1+\frac{p}{N+p}}\\
    \leq &\,c \left(\frac{\varepsilon}{\theta}\right)^\frac{p}{N+p}2^{j[2p-(p-2-q)\frac{p}{N+p}]}\rho^{-p}\varepsilon^{-1}|A_j|^{1+\frac{p}{N+p}}\\
    \leq &\,c \varepsilon^{-\frac{N}{N+p}}\theta^{-\frac{p}{N+p}} b^j\rho^{-p}\varepsilon^{-1}|A_j|^{1+\frac{p}{N+p}}
\end{align*}
with $b=2^{[2p-(p-2-q)\frac{p}{N+p}]}$. Since $|Q_j|=c(N)\theta \rho^{N+p}$, we get
\[
    Y_{j+1}\leq c \varepsilon^{-\frac{N}{N+p}}b^jY_j^{1+\frac{p}{N+p}}.
\]
Setting $\delta=p/(N+p)$ we see that Lemma \ref{Lm:fast_convergence} guarantees  $Y_j\to 0$ provided
\[
    \frac{|Q_{\rho,\theta \rho^p}\cap \{u>\mu_+/2\}|}{|Q_{\rho,\theta\rho^p}|}=Y_o\leq\left(c \varepsilon^{-\frac{N}{N+p}}\right)^{-\frac{1}{\delta}}b^{-\frac{1}{\delta}}=:\varepsilon^\frac{N}{p}\nu_o,
\]
where $\nu_o$ depends only on $p,q,N$. Since $|Q_j|$ is bounded from above, this means that $|A_j|\to 0$. Furthermore, since
\[
    Q_{\rho/2,\theta(\rho/2)^p}\cap \{u>\frac{3}{4}\mu_+\}\subset A_j
\]
for every $j$, the measure of the set on the left hand side must be zero.
\end{proof}

\begin{remark}\label{Rmk:A1:1}
{\normalfont If we let $\rho \equiv R_o$, since $Q_\rho(\theta) \equiv Q_1$ and $Q_1 \subseteq Q_o$, 
it is apparent that \eqref{Eq:A1:4} is satisfied, and therefore, the results of Lemma~\ref{Lm:A1:1} apply to $Q_1$.}
\end{remark}
\begin{remark}\label{Rmk:A1:2}
{\normalfont It is straightforward to check from the proof that Lemma~\ref{Lm:A1:1} holds for any $q>0$, and any $p>1$, and 
not only for the indicated ranges of $p$ and $q$.}
\end{remark}
\noindent For simplicity, from here on we write $\mu_{+}$ and $\mu_{-}$, instead of $\mu_o^{+}$and $\mu_o^{-}$.
Moreover, as already pointed out in Remark~\ref{Rmk:A1:1}, we let $\rho \equiv R_o$, so that Lemma~\ref{Lm:A1:1} can be applied.
\subsection{The Reduction of the Oscillation}
We are now ready to prove the reduction of oscillation in the case $\mu_{-}=0$, which implies that $\mu_{+}=\omega\equiv\om_o$.
If \eqref{Eq:A1:5} applies, by Lemma~\ref{Lm:A1:1} we have a reduction of oscillation from above, i.e.
\[
\essosc_{(x_o,t_o)+Q_{\frac\rho2}(\theta)} u \le \omega_1, \qquad \omega_1\df=\frac{3}{4} \omega_o.
\]
Suppose now that \eqref{Eq:A1:5} fails, i.e. that
\begin{equation}\label{Eq:A1:7}
\left|\left\{u>\frac{\mu_{+}}2\right\} \cap(x_o,t_o)+Q_\rho(\theta)\right|>\nu_o \delta^{N/p}\left|Q_\rho(\theta)\right|.
\end{equation}
This implies that there is a set $A\subset(t_o-\theta\rho^p,t_o]$ of positive measure, such that
\begin{equation}\label{Eq:A1:8}
\left|\left\{x \in K_\rho(x_o):\, u(x,\tau)>\frac{\mu_+}{2}\right\}\right|>\nu_o \delta^{N/p}\left|K_\rho\right|
\end{equation} 
for all $\tau\in A$. Indeed, if \eqref{Eq:A1:8} failed for a.e. $\tau\in(t_o-\theta\rho^p,t_o]$, \eqref{Eq:A1:7}
would be violated as well.

Let us now recall the following $L^1$-Harnack inequality.
\begin{proposition}[Proposition 7.1 of \cite{BDGLS}]\label{Prop:A1:1}
Assume that $0<p-1<q$. There exists a positive constant $\gamma$ depending only on the data $\left\{q, p, N, C_o, C_1\right\}$ such that whenever $u$ is a 
non-negative, local, bounded, weak solution to \eqref{eq:structure}--\eqref{Def:structure} in $\Omega_T$, then for every cylinder 
$\dsty K_{2\rho}(y)\times[s,\tau] \subset\subset \Omega_T$, we have
$$
\sup _{t \in[s,\tau]} \int_{K_\rho(y) \times\{t\}} u^q\, dx \le \gamma \inf_{t \in[s,\tau]} \int_{K_{2\rho}(y) \times\{t\}} u^q\, dx 
 +\gamma\left(\frac{\tau-s}{\rho^\lambda}\right)^{\frac{q}{q+1-p}},
$$
where $\lambda=\frac{N}{q}(p-q-1)+p$.
\end{proposition}
Provided that $(x_o,t_o)+Q_{2\rho}(\theta) \subset \Omega_T$, for $\tau\in A$ we have
$$
\begin{aligned}
\int_{K_\rho(x_o)} u^q(x,\tau)\,dx&\ge\left|\left\{x \in K_\rho(x_o):\, u(x,\tau)>\frac{\mu_+}{2}\right\}\right|\left(\frac{\mu_{+}}{2}\right)^q\\ 
&>\left(\frac{\mu_+}{2}\right)^q \nu_o \delta^{N/p}(2 \rho)^N.
\end{aligned}
$$
Hence
\begin{align*}
\left(\frac{\mu_+}{2}\right)^q \nu_o \delta^{N/p}(2 \rho)^N\le&\gamma \inf_{t \in[t_o-\theta\rho^p,t_o]} \int_{K_{2\rho}(x_o)\times\{t\}} u(x,t)^q\, dx\\
&+\gamma\left[\frac\theta{\rho^{\frac Nq(p-q-1)}}\right]^{\frac q{q+1-p}}.
\end{align*}
Notice that
$$
\left[\frac{\theta}{\rho^{\frac{N}{q}(p-q-1)}}\right]^{\frac{q}{q+1-p}}=\frac{[\delta \mu_+^{q+1-p}]^{\frac q{q+1-p}}}{\rho^{-\frac{N}{q} \frac{q+1-p}{q+1-p} q}}
=\delta^{\frac{q}{q+1-p}} \mu_{+}^q \rho^N
$$
This yields
\begin{align*}
\gamma \inf_{t \in[t_o-\theta\rho^p,t_o]} \int_{K_{2\rho}(x_o)\times\{t\}} u(x,t)^q\, dx&\ge\mu_+^q\rho^N\left[\nu_o \delta^{N/p} 2^{N-q}-\gamma\delta^{\frac{q}{q+1-p}} \right]\\
&\ge\mu_+^q\rho^N\delta^{N/p}\left[\nu_o 2^{N-q}-\gamma\delta^{\frac{q}{q+1-p}-\frac Np} \right].
\end{align*}
Moreover, we have $\displaystyle \delta^{\frac{q}{q+1-p}-\frac{N}{p}}=\delta^{\frac{qp-N(q+1-p)}{p(q+1-p)}}$.
We wonder when $qp-N(q+1-p)>0$, which is the same as inquiring when $q(N-p)-N(p-1)<0$. If $p < N$ we have $q(N-p)<N(p-1)$, which yields $q<\frac{N(p-1)}{N-p}$. 
If $p \ge N$ then $q(p-N)>N(1-p)$ is always true.
Hence, we can say that the condition
$$
p-1<q<\frac{N(p-1)}{(N-p)_+}
$$
covers the range in which the exponent we are interested in is positive. Hence, if we choose $\delta$ such that
$$
\nu_o 2^{N-q}-\gamma \delta^{\frac{N(p-1)-q(N-p)}{p(q+1-p)}}=\frac12\nu_o 2^{N-q} \quad\Rightarrow\quad
\delta^{\frac{N(p-1)-q(N-p)}{p(q+1-p)}}=\frac{\nu_o 2^{N-q}}{2\gamma},
$$
we obtain
$$
C(\text{data}) \mu_{+}^q \rho^N \leq \inf_{t\in[t_o-\theta \rho^p,t_o]}  \int_{K_{2\rho}(x_o)} u^q(x,t)\,dx
$$
where $C(\text{data})=C(q,p,N,C_o,C_1)$.
Now, take $\eta\in(0,1)$ and remark that
\[
\begin{aligned}
\int_{K_{2\rho}(x_o)} u^q(x,t)\,dx=&\int_{K_{2\rho}(x_o)\cap[u\ge\eta\mu_+]} u^q(x,t)\,dx\\
&+\int_{K_{2\rho}(x_o)\cap[u<\eta\mu_+]} u^q(x,t)\,dx \\
& \le \mu_{+}^q\left|K_{2\rho}(x_o)\cap\left[u(\cdot,t) \ge \eta \mu_{+}\right]\left|+\eta^q \mu_{+}^q\right|K_{2\rho}(x_o)\right|.
\end{aligned}
\]
Hence, we have
\[
\begin{aligned}
&C \mu_+^q \rho^N \le \mu_{+}^q \left|K_{2\rho}(x_o) \cap[u(\cdot,t) \ge \eta \mu_+]\right|+\eta^q\mu_+^q \left|K_{2\rho}(x_o)\right|\\
&C_1\left|K_{2\rho}(x_o)\right| \le \left|K_{2\rho}(x_o)\cap\left[u(\cdot,t) \ge \eta \mu_+\right]\right|+\eta^q\left|K_{2\rho}(x_o)\right|\\
&\left(C_1-\eta^q\right)\left|K_{2\rho}(x_o)\right| \le \left|\left[u(\cdot,t) \ge \eta \mu_+\right] \cap K_{2\rho}(x_o)\right|.\\
\end{aligned}
\]
If we choose $\eta$ such that
$$
C_1 -\eta^q=\frac{1}{2} C_1\quad\Rightarrow\quad \frac{1}{2} C_1=\eta^q\quad\Rightarrow\quad \eta=\left(\frac{1}{2} C_1\right)^{1/q},
$$
we conclude that $\forall t \in\left(t_o-\theta\rho^p,t_o\right]$
\begin{equation}\label{Eq:A1:9}
\left|\left[u(\cdot,t) \ge \eta \mu_+\right] \cap K_{2\rho}(x_o)\right| \ge \alpha\left|K_{2\rho}(x_o)\right|
\end{equation}
for a constant $\alpha$ that depends only on the data. 

We rely on the expansion of positivity of Proposition~\ref{Prop:exp_pos}; to avoid any kind of confusion, here we denote with $\bar\delta$ and $\eta_*$ the quantities $\delta$ and $\eta$ of the proposition.
Provided we take
$$
M=\min \left\{\eta \mu_+,\left(\frac{\delta}{2^p}\right)^{\frac1{q+1-p}} \mu_+\right\},
$$
we have that \eqref{Eq:A1:9} still holds with $\eta\mu_+$ replaced by $M$, and also that
$$
K_{32\rho}(x_o) \times(t_o-\bar{\delta} M^{q+1-p}(2\rho)^p,t_o] \subset
K_{32 \rho}(x_o) \times(t_o-\theta\rho^p, t_o] \subset \Omega_T,
$$
where, as just mentioned above, $\bar\delta\in(0,1)$ is the constant which appears in Proposition~\ref{Prop:exp_pos}.
Hence, we may apply Proposition~\ref{Prop:exp_pos} with
$$
s=t_o-\bar{\delta} M^{q+1-p}(2 \rho)^p
$$
and $\rho$ replaced by $2\rho$ to conclude that there is a $\xi\in(0,1)$ and $\tilde{\delta}<\bar\delta$ depending only on the data, such that
$$
u>\xi\mu_{+} \quad \text { in } K_{4 \rho}(x_o) \times(t_o-\tilde{\delta} \mu_{+}^{q+1-p} \rho^p,t_o],
$$
which yields the reduction of oscillation from below.
Combining the previous reasoning and Lemma~\ref{Lm:A1:1}, we have the following.
\begin{lemma}\label{Lm:A1:3} 
There exist constants $\delta, \gamma, \eta \in(0,1)$ that depend only on the data, such that for any weak solution $u$ and number $\mu_+$ satisfying
$$
K_{32\rho}(x_o) \times(t_o-\delta \mu_{+}^{q+1-p} \rho^p, t_o] \subset \Omega_T
$$
and
$$
u \le \mu_{+} \quad\text{ in } \quad K_\rho(x_o) \times(t_o-\delta \mu_{+}^{q+1-p} \rho^p, t_o],
$$
we have
\begin{equation}\label{Eq:A1:10}
\essosc_{K_{\frac{\rho}{2}}(x_o) \times(t_o-\gamma \delta \mu_+^{q+1-p} \rho^p, t_o]} u \le \eta \mu_{+}.
\end{equation}
Furthermore, if we let for simplicity $Q_1\df=K_{\frac{\rho}{2}}(x_o) \times(t_o-\gamma \delta \mu_+^{q+1-p} \rho^p, t_o]$, one of the following conditions must hold:
\begin{align}
either\qquad&\esssup_{Q_1} u\le\frac{1+\eta}{2} \mu_+\label{Eq:A1:11}\\
or\qquad&\essinf_{Q_1} u\ge \frac{1-\eta}{2} \mu_+\label{Eq:A1:12}.
\end{align}
\end{lemma}
\begin{proof}
The quantity $\eta$ in \eqref{Eq:A1:10} is given by $\dsty\max \left\{(1-\xi),\frac{3}{4}\right\}$,
according to which of the two alternatives holds. The quantity $\gamma$ (which measures the width of the time gap) is given by
$\dsty\min \left\{2^{-p}, \frac{\tilde\delta}{\delta}\right\}$. If \eqref{Eq:A1:11} holds true, there is nothing to prove; if \eqref{Eq:A1:11} does not hold true, we have
$$
\begin{aligned}
\esssup_{Q_1} u&-\essinf_{Q_1} u\le\eta\mu_+,\\
\essinf_{Q_1} u&\ge \esssup _{Q_1} u - \eta \mu_+ \ge \frac{1+\eta}{2} \mu_+-\eta \mu_+=\frac{1-\eta}{2} \mu_{+}, 
\end{aligned}
$$
which is precisely \eqref{Eq:A1:12}.
\end{proof}
\subsection{The Reduction of the Oscillation Concluded}
We define
\[
C:=2\max\{2,\gamma^{-\frac1p}\},\qquad\mu_+^j:=\left(\frac{1+\eta}2\right)^j \mu_+,\qquad\rho_j:=\frac\rho{C^j},
\]
where $\gamma$ is the constant claimed in Lemma~\ref{Lm:A1:3}. Once $\rho_j$ and $\mu_+^j$ are given in this way, it is apparent that 
\[
Q_{\rho_1,\delta(\mu_+^1)^{q+1-p}\rho_1^p}\subseteq Q_{\frac\rho2,\gamma\delta(\mu_+)^{q+1-p}\rho^p}
\]
and 
\[
\essosc_{Q_{\rho_1,\delta(\mu_+)^{q+1-p}\rho_1^p}}u\le\eta\mu_+.
\]
Moreover, if \eqref{Eq:A1:11} occurs, we also have that $\dsty u\le\mu_+^1$ in $Q_{\rho_1,\delta(\mu_+^1)^{q+1-p}\rho_1^p}$. We may then 
apply Lemma~\ref{Lm:A1:3} once more, and conclude not only that 
\[
\essosc_{Q_{\rho_2,\delta(\mu_+^1)^{q+1-p}\rho_2^p}}u\le\eta\mu_+^1,
\]
but also that either \eqref{Eq:A1:11} or \eqref{Eq:A1:12} holds with $\rho$ replaced by $\rho_1$, and $\mu_+$ replaced by $\mu_+^1$. If \eqref{Eq:A1:11}
holds true, we can repeat the whole procedure one more time. Hence, provided that \eqref{Eq:A1:11} continue to hold, we end up with 
\begin{align}
\essosc_{Q_{\rho_j,\delta(\mu_+^{j-1})^{q+1-p}\rho_j^p}}u&\le\eta\mu_+^{j-1},\label{Eq:A1:13}\\
\esssup_{Q_{\rho_{j-1},\delta(\mu_+^{j-1})^{q+1-p}\rho_{j-1}^p}}u&\le\mu_+^{j-1}.\label{Eq:A1:14}
\end{align}
Either this procedure continues indefinitely $\forall\,j\in\Nn$, or there exists $k\in\Nn$ such that \eqref{Eq:A1:13} holds for all $j=1,2,\dots,k$, and 
\begin{equation}\label{Eq:A1:15}
\essinf_{Q_{\frac{\rho_{k-1}}2,\gamma\delta(\mu_+^{k-1})^{q+1-p}\rho_{k-1}^p}}u\ge\left(\frac{1-\eta}2\right)\mu_+^{k-1}=\left(\frac{1-\eta}{1+\eta}\right)\mu_+^k.
\end{equation}
Let us see the consequences of \eqref{Eq:A1:15}. We remark that
\begin{align*}
&\frac{\rho_{k-1}}2=\frac C2\rho_k\ge 2\rho_k\\
&\gamma\delta\left(\mu_+^{k-1}\right)^{q+1-p}\rho_{k-1}^p=\gamma C^p\left(\frac2{1+\eta}\right)^{q+1-p}\delta\left(\mu_+^k\right)^{q+1-p}\rho_k^p
\ge 2^p \left(\frac2{1+\eta}\right)^{q+1-p}\delta\left(\mu_+^k\right)^{q+1-p}\rho_k^p\\
&\hskip3.truecm>\delta \left(\mu_+^k\right)^{q+1-p}(2\rho_k)^p.
\end{align*}
Hence, it follows from \eqref{Eq:A1:14} and \eqref{Eq:A1:15} with $j=k$ that 
\[
\left(\frac{1-\eta}{1+\eta}\right)\mu_+^k\le u\le\frac2{1+\eta}\mu_+^k\qquad\text{ in }\,Q_{2\rho_k,\delta(\mu_+^k)^{q+1-p}(2\rho_k)^p}.
\]
This is exactly the assumption of \cite[Lemma~3.7]{VV_2022} and it shows that the function 
\[
v(x,t):=\left(\mu_+^{k-1}\right)^{-1}u(x,\left(\mu_+^k\right)^{q+1-p}t),\qquad (x,t)\,\text{ in }\,Q_{2\rho_k,\delta(2\rho_k)^p}(x_o,0),
\]
solves a parabolic equation of parabolic $p$-laplacian-type, for which it is known that there exist constants $\epsilon$,
$\delta\in(0,1)$, and $A\ge1$, 
depending only on $p,q,N$ such that, setting
\begin{equation*}
\om_n = \delta\om_{n-1}, \quad
\theta_n = \left ( \frac{\om_n}A \right )^{q+1-p}, \quad
\rho_n = \epsilon \rho_{n-1} \quad \textrm{and} \quad
Q_n = Q_{\rho_n}^- (\theta_n),
\end{equation*}
for all non-negative integers $n$, there holds $Q_{n+1}\subset Q_n$ and
\begin{equation}\label{Eq:A:1:16}
\osc\limits_{Q_n} u \le \om_n.
\end{equation}
Therefore, either \eqref{Eq:A1:13} holds for any $j\in\Nn$, which we can rephrase as in \eqref{Eq:A:1:16}, or there exists $k\in\Nn$ for which \eqref{Eq:A1:14} fails, and in this case again we have \eqref{Eq:A:1:16} for every $n>k$. In any case, it is apparent that we have precisely the content of Lemma~\ref{Lm:Geom-Decay}, which is therefore proved.

From here on, working exactly as in \cite[Lemma~8.3 and Theorem~8.4]{VV_2022}, one concludes that $u$ is locally H\"older continuous in $\Omega_T$, and that the H\"older exponent depends on $p$, $q$, $N$, $C_o$, $C_1$. We refrain from going into details here, since the arguments are essentially the same and could be reproduced verbatim.
\section{H\"older continuity at the Boundary}\label{A:S:Bd}
As stated in section~\ref{S:A:Intro}, we assume that $u$ vanishes at the lateral boundary $S_T$ and that $\partial\Om$ satisfies a geometric density condition.
We extend $u$ by zero outside $\Om_T$; it can then be proved that $u$ is a weak sub-solution. We refrain from going into details about this property, since it is a standard argument. In the following, we continue to denote by $u$ such an extension.

We repeat and adapt at the boundary the same construction of initial cylinders, as in Section~\ref{SS:A:B2}.
Namely, fix $(x_o,t_o)\in S_T$, and for $\delta \in(0,1)$ as in Lemma~\ref{Lm:A1:3}, and $0<R_o<R_*$ sufficiently small such that $(t_o-\delta R_o^p, t_o]\subset(0,T]$,
define
\begin{equation}\label{Eq:C1:3}
Q_o^*\df=K_{R_o}(x_o) \times(t_o-\delta R_o^p, t_o], 
\end{equation}
and let
$$
\tilde{\mu}_o^+= \esssup _{Q_o^{*}} u,\quad\tilde{\omega}_o=\essosc_{Q_o^*} u=\tilde{\mu}_o^{+}.
$$
Recall that by construction $\quad\tilde{\mu}_o^-= \essinf_{Q_o^{*}} u=0$. Since $u$ is locally bounded, without loss of generality we may assume that $\tilde{\omega}_o \le 1$, so that
$$
Q_o\df=K_{R_o}(x_o) \times(t_o- \tilde{\omega}_o^{q+1-p} \delta R_o^p, t_o] \subseteq Q_o^* \subset \Omega_T.
$$
We let $\dsty\mu_o^{+}=\esssup_{Q_o} u$ (again $\dsty\mu_o^{-}=\essinf_{Q_o} u=0$). Since $\mu_o^{+} \le \tilde{\mu}_o^{+}$, it is apparent that if we let
$$
Q_1\df=K_{R_o}(x_o) \times(t_o-\delta\left(\mu_o^{+}\right)^{q+1-p} R_o^p, t_o],
$$
then $\dsty Q_1 \subseteq Q_o$. 
\vskip.2truecm
The function $x\mapsto (u(x,t)-k)_+$ vanishes outside $(\partial\Om\cap B_\rho)$ for all $t\in(t_o-\delta\left(\mu_o^{+}\right)^{q+1-p}\rho^p,t_o]$ with $0<\rho\le R_o$, and $\partial\Om$ satisfies the geometric density condition \eqref{Eq:11:8:4}; this yields the following.
\begin{lemma}\label{Lm:C:1} 
Let $u$ be a non-negative, local, bounded, weak sub-solution to \eqref{eq:structure}--\eqref{Def:structure} in $\Omega_T$ for 
$0<p-1<q<\frac{N(p-1)}{(N-p)_{+}}$. Set
$\dsty\theta=\delta(\mu_o^{+})^{q+1-p}$
for the same $\delta \in(0,1)$ as in Lemma~\ref{Lm:A1:3}. Assume that $(x_o,t_o)\in S_T$, $\dsty [(x_o,t_o)+Q_\rho(\theta)]\cap\Om_T \subset Q_1$, and
\begin{equation}\label{Eq:C1:4}
\esssup_{(x_o,t_o)+Q_\rho(\theta)} u \le \mu_o^{+}.
\end{equation}
Then, $\forall\,t\in(t_o-\delta\left(\mu_o^{+}\right)^{q+1-p}\rho^p,t_o]$, and $\forall\,s\ge1$
\[
\left|\left[u(\cdot,t)\ge\mu_o^+-\frac\om{2^s}\right]\cap B_\rho(x_o)\right|\le\alpha|B_\rho|,
\]
where $\alpha=1-\beta$.
\end{lemma}
We can now repeat almost verbatim the proof of \cite[Lemma~6.3]{BDGLS}, properly adjusted to our context, and conclude that for any $\nu\in(0,1)$ there exists $s_*\in\Nn$, $s_*\ge2$, $s_*$ independent of $\mu_o^+$ and $\rho$, such that
\[
\left|\left[u\ge\mu_o^+-\frac\om{2^{s_*}}\right]\cap Q_\rho\left(\frac12\theta\right)\right|\le\nu\,\left|Q_\rho\left(\frac12\theta\right)\right|.
\]
In Lemma~\ref{Lm:A1:1}, $u>\frac{\mu_o^+}2$ in \eqref{Eq:A1:5} implies $u<\frac34\mu_o^+$ in \eqref{Eq:A1:6}; however, a careful check of the proof shows that $u>\mu_o^+-\frac{\om}{2^{s_*}}$ in \eqref{Eq:A1:5} implies $u<\mu_o^+-\frac{\om}{2^{s_*+1}}$ in \eqref{Eq:A1:6} (indeed, it is just a matter of simple computations). Hence, a straightforward application of Lemma~\ref{Lm:A1:1} with $\delta$ replaced by $\frac12\delta$, allows to conclude that 
\[
u\le\mu_o^+\left(1-\frac1{2^{s_*+1}}\right)\qquad\text{ a.e. in }\,\left[(x_o,t_o)+Q_{\frac12\rho}\left(\frac12\theta\right)\right]\cap \Om_T.
\]
This plays the role of \eqref{Eq:A1:10} as far as the reduction of the oscillation at the boundary is concerned. Moreover, since $(x_o,t_o)\in S_T$, the second alternative never occurs, and we can indefinitely apply the iteration and conclude.

The gist of the whole argument is that once more we end up with a geometric decay of $u$, exactly as given in Lemma~\ref{Lm:Geom-Decay}, with the only difference that for $(x_o,t_o)\in S_T$ \eqref{Eq:geom-decay} is substituted by \eqref{Eq:geom-decay-bd}, as already pointed out at the beginning.
\section{Construction of Signed Solutions by Odd Reflection}\label{S:A:Ext}
We prove that we can extend a positive solution in the upper-half space to a signed solution to a similar equation in the whole space by odd reflection.

\begin{lemma}\label{Lm:odd_reflection}
    Let 
    \[
    Q^+ = \{(x,t) : |x_i|<1, 0<x_N<2, t_1<t<t_2\}
    \]
    and
    \[
    Q^- = \{(x,t) : |x_i|<1, -2<x_N<0, t_1<t<t_2\}.
    \]
    Let $u$ solve \eqref{eq:structure}--\eqref{Def:structure} in $Q^+$ with $u=0$ on $\partial Q^+\cap \{x_N=0\}$. Define
    \[
    Q' = \{(x,t) : |x_i|<1, -2<x_N<2, t_1<t<t_2\}
    \]
    and $\tilde u, \tilde A$ as
    \begin{align*}
        \tilde u(x', x_N, t)=\begin{cases}
            u(x', x_N, t), &\mbox{ if } x_N\geq 0,\\
            -u(x', -x_N, t), &\mbox{ if } x_N< 0,
        \end{cases}
    \end{align*}
    \begin{align*}
        \tilde A_i(x', x_N, t)=\begin{cases}
            A_i(x', x_N, t), &\mbox{ if } x_N\geq 0,\\
            -A_i(x', -x_N, t), &\mbox{ if } x_N< 0,
        \end{cases}, \, i= 1, ..., N-1,
    \end{align*}
    and
    \begin{align*}
        \tilde A_N(x', x_N, t)=\begin{cases}
            A_N(x', x_N, t), &\mbox{ if } x_N\geq 0,\\
            A_N(x', -x_N, t), &\mbox{ if } x_N< 0,
        \end{cases}
    \end{align*}
    where $\tilde A(x', x_N, t):= \tilde A(x', x_N, t, u(x', x_N, t), Du(x', x_N, t))$.

    Then $\tilde u$ is a weak solution to \eqref{Eq:A1:1bis} in $Q'$ with $A$ replaced by $\tilde A$.
    
\end{lemma}
\begin{proof}
    We follow the proof of \cite[Lemma 2.7]{lundstrom2009}. Clearly we still have
    \[
    \tilde u\in C_{\loc}(t_1,t_2;L_{\loc}^{q+1}(K'))\cap L_{\loc}^p(t_1,t_2;W_{\loc}^{1,p}(K')),
    \]
    where
    \[
    K^+ = \{x : |x_i|<1, 0<x_N<2\}
    \]
    and $K^-, K'$ are defined similarly. Let $\theta\in C^\infty_o(Q')$. We aim at establishing that
\begin{align}\label{eq:odd_thesis}
    \iint_{Q'} -|\tilde u|^{q-1}\tilde u\partial_t \theta +\tilde A(x', x_N, t)\cdot D\theta \, dxdt =0.
\end{align}
    For this purpose, write $\theta = \phi + \psi$, where
    \[
        2\phi(x)=\theta(x',x_N)+\theta(x',-x_N)
    \]
    is even in $x_N$, and thus
    \[
        2\psi(x)=\theta(x',x_N)-\theta(x',-x_N)
    \]
    is odd in $x_N$. Furthermore $\psi=0$ in $Q'$ when $x_N=0$ and so it belongs to 
    \[
    \psi\in W_{\loc}^{1,q+1}(t_1,t_2;L^{q+1}(K+))\cap L_{\loc}^p(t_1,t_2;W^{1,p}_o(K+)).
    \]
    From the equation satisfied by $u$ we get immediately
    \[
        \iint_{Q^+} -|\tilde u|^{q-1}\tilde u\partial_t \psi +\tilde A(x', x_N, t)\cdot D\psi \, dxdt =0.
    \]
    Now we check this condition in $Q^-$.
\begin{align*}
    &\iint_{Q^-} -|\tilde u|^{q-1}\tilde u\partial_t \psi +\tilde A(x', x_N, t)\cdot D\psi \, dxdt \\
        =&\,\iint_{Q^-} -|- u(x',-x_N,t)|^{q-1}( -u(x',-x_N,t))\partial_t \psi \\
        &\,+\sum_{i=1}^{N-1} -A_i(x',-x_N,t,\tilde u(x',-x_N,t),D\tilde u(x',-x_N,t))\partial_i \psi\\
        &\,+ A_N(x',-x_N,t,\tilde u(x',-x_N,t),D\tilde u(x',-x_N,t)) \partial_N \psi\, dxdt
\end{align*}
Note that in $Q^-$, $\tilde u(x',-x_N,t)=u(x',-x_N,t)$ and $D\tilde u(x',-x_N,t)=D u(x',-x_N,t)$. Now we make the change of variables $x_N \to -x_N$ and get $Q^-\to Q^+$, $dx\to -dx$, and use that $\psi$ is odd in $x_N$, $\partial_i\psi$ is odd in $x_N$ for $i=1, ..., N-1$ and even for $i=N$.
\begin{align*}
    &\iint_{Q^-} -|\tilde u|^{q-1}\tilde u\partial_t \psi +\tilde A(x', x_N, t)\cdot D\psi \, dxdt \\
    =&\,-\iint_{Q^+} | u(x',x_N,t)|^{q-1}( u(x',x_N,t))\partial_t (-\psi(x',x_N,t)) \\
        &\,+\sum_{i=1}^{N-1} A_i(x',x_N,t, u(x',x_N,t),D u(x',x_N,t))\partial_i \psi\\
        &\,+ A_N(x',x_N,t, u(x',x_N,t),D u(x',x_N,t)) \partial_N \psi\, dxdt \\
    =&\, \iint_{Q^+} |u|^{q-1}u \partial_t \psi-A(x,t,u,Du)\cdot D\psi\,dxdt=0.
\end{align*}
Similarly we get the same result for $\phi$, noting that $\partial_i \phi$ is even in $x_N$ for $i=1, ..., N-1$ and $\partial_N \phi$ is odd. Combining both we get \eqref{eq:odd_thesis}, as intended. 
\end{proof}
\section{Oscillation Decay for Signed Solutions}\label{A:S:Geom-signed}
Suppose $u$ is a signed solution built by odd reflection along a flat boundary of a non-negative solution originally defined in $\Om_T$ (for simplicity, we do not distinguish between the original solution, and its extension).

If the starting cylinder $(x_o,t_o)+Q_{\rho_o}^-(\theta_o)\subset\subset\Om_T$, the geometric decay of $u$ is described by Lemma~\ref{Lm:Geom-Decay}.  By construction, we have the same behavior for the non-positive $u$ in $(x_o,t_o)+Q_{\rho_o}^-(\theta_o)\subset\subset(\Om_T)^c$.

When $(x_o,t_o)\in S_T$, the controlled reduction of the oscillation of $u$ is provided by the arguments of Section~\ref{A:S:Bd}. Taking the worst parameters, we eventually achieve the wanted decay for the signed solution we are interested in.


\begin{thebibliography}{99}
\bibitem{adamowicz} T. Adamowicz and N.L.P. Lundstr\"om, The boundary Harnack inequality for variable exponent $p$-Laplacian, Carleson estimates, barrier functions and $p(\cdot)$-harmonic measures, \textit{Ann. Mat. Pura Appl. (4), \textbf{195}, (2016), 623--658.}

\bibitem{Aikawa} H. Aikawa, T. Kilpel\"{a}inen, N. Shanmugalingam and X. Zhong, {Boundary {H}arnack principle for {$p$}-harmonic functions in
smooth {E}uclidean domains}, \textit{Potential Anal.,
\textbf{26}(3), (2007), 281--301.}

\bibitem{ACS1} I. Athanasopoulos, L. Caffarelli and S. Salsa, {Caloric functions in Lipschitz domains and the regularity of
solutions to phase transition problems}, \textit{Ann. of Math. (2), \textbf{143}(3), (1996), 413--434.}

\bibitem{ACS2} I. Athanasopoulos, L. Caffarelli and S. Salsa, {Regularity of the free boundary in parabolic phase-transition
problems}, \textit{Acta Math. \textbf{176}(2), (1996), 245--282.}

\bibitem{ACS3} I. Athanasopoulos, L. Caffarelli and S. Salsa, {Phase transition problems of parabolic type: flat free boundaries are smooth}, \textit{Comm. Pure Appl. Math., \textbf{51}(1), (1998), 77--112.}

\bibitem{AGS} B. Avelin, U. Gianazza and S. Salsa, 
Boundary estimates for certain degenerate and
singular parabolic equations, \textit{J. Eur. Math. Soc. \textbf{18}, (2016), 381--424.}


\bibitem{bogelein2021} V. Bögelein, F. Duzaar and N. Liao, On the Hölder regularity of signed solutions to a doubly nonlinear equation, \textit{J. Funct. Anal., \textbf{281}(9), (2021), article no. 109173, 58 pp.}

\bibitem{bogelein2022} V. Bögelein, F. Duzaar, N. Liao and L. Schätzler, On the Hölder regularity of signed solutions to a doubly nonlinear equation. Part II, \textit{Rev. Mat. Iberoam., \textbf{39}(3), (2022), 1005--1037.}

\bibitem{BS} V. Bögelein and M. Strunk, A comparison principle for doubly nonlinear parabolic partial differential equations, \textit{Ann. Mat. Pura Appl. (4), (2023), https://doi.org/10.1007/s10231-023-01381-4.}

\bibitem{BDGLS} V. Bögelein, F. Duzaar, U. Gianazza, N. Liao and C. Scheven, H\"older Continuity of the Gradient of Solutions to Doubly Non-Linear Parabolic Equations, \textit{arXiv preprint, (2023), 1--142, doi:10.48550/arxiv.2305.08539.}


\bibitem{Caf87}
L. Caffarelli, 
{A Harnack inequality approach to the regularity of free boundaries. I. Lipschitz free boundaries are $C^{1,\alpha}$},
\textit{Rev. Mat. Iberoamer. \textbf{3}, (1987), 139--162.}

\bibitem{Caf89}
L. Caffarelli, 
{A Harnack inequality approach to the regularity of free boundaries. II. Flat free boundaries are Lipschitz},
\textit{Comm. Pure Appl. Math. \textbf{42}, (1989), 55--78.}


\bibitem{CS05} L. Caffarelli and S. Salsa, \textit{A geometric approach to free boundary problems}, Graduate Studies in Mathematics, 68,
{American Mathematical Society, Providence, RI}, 2005.

\bibitem{De11}
D. De Silva, 
{Free boundary regularity for a problem with right hand side.},
\textit{Interfaces Free Bound., \textbf{13}(2), 2011, 223--238.}

\bibitem{DFS21}
D. De Silva, N. Forcillo and O. Savin, 
Perturbative estimates for the one-phase Stefan problem,
\textit{Calc. Var., \textbf{60}(6), 38, (2021), 1--38.}

\bibitem{DiBen_1993}
E.  DiBenedetto, 
\textit{Degenerate Parabolic Equations},
Universitext, Springer, New York, 1993.

\bibitem{DGV}
 E. DiBenedetto, U. Gianazza and V. Vespri, 
\textit{Harnack’s inequality for degenerate and singular parabolic},
Springer Monographs in Mathematics, Springer-Verlag, New York, 2012.

\bibitem{Fer06}
 F. Ferrari, 
Two-phase problems for a class of fully nonlinear elliptic operators. Lipschitz free
boundaries are $C^{1,\gamma}$,
\textit{Amer. J. Math. \textbf{128}, (2006), 541--571.}

\bibitem{FS07}
 F. Ferrari and S. Salsa, 
Regularity of the free boundary in two-phase problems for linear elliptic operators,
\textit{Adv. Math. \textbf{214}, (2007), 288--322.}

\bibitem{FS10}
 F. Ferrari and S. Salsa, 
Regularity of the solutions for parabolic two-phase free boundary problems,
\textit{Comm. Partial Differential Equations \textbf{35}, (2010),  1095--1129.}

\bibitem{kuusi2014} T. Kuusi, G. Mingione and K. Nyström, A boundary Harnack inequality for singular equations of $p$-parabolic type, \textit{Proc. Amer. Math. Soc., \textbf{142}(8), (2014), 2705--2719.}

\bibitem{kuusi2020} T. Kuusi, M. Misawa and K. Nakamura, Regularity Estimates for the $p$-Sobolev Flow, {\it J. Geom. Anal. \textbf{30}, (2020), 1918--1964.}

\bibitem{liao2021} N. Liao and L. Schätzler, On the Hölder regularity of signed solutions to a doubly nonlinear equation. Part III, \textit{Int. Math. Res. Not., \textbf{3}, (2022), 2376--2400.}

\bibitem{lundstrom2009} N.L.P. Lundstr\"om and K. Nystr\"om, On a two-phase free boundary condition for p-harmonic
measures, \textit{ manuscripta math., \textbf{129}, (2009), 231--249.}

\bibitem{moring2023supercaloric} K. Moring and C. Scheven, Supercaloric functions for the porous medium equation in the fast diffusion case, \textit{arXiv preprint, (2023), 1--34, doi:10.48550/arXiv.2306.07155.}

\bibitem{SV_1994}
G. Savar\'{e} and V. Vespri, 
The asymptotic profile of solutions of a class of doubly
nonlinear equations,
\textit{Nonlinear Anal., \textbf{22}(12), (1994), 1553--1565.}

\bibitem{VV_2022} V. Vespri and M. Vestberg,
An extensive study of the regularity of solutions to doubly singular equations,
\textit{Adv. Calc. Var., \textbf{15}(3), (2022), 435--473.}

\end{thebibliography}
\end{document}